\documentclass[11pt]{amsart}
\usepackage[margin=1in]{geometry}
\usepackage{amsmath,amssymb,amsfonts,amsthm,mathtools}
\numberwithin{equation}{section}
\usepackage[colorlinks=true,linkcolor=blue,citecolor=blue,urlcolor=blue]{hyperref}
\theoremstyle{plain}
\newtheorem{theorem}{Theorem}[section]
\newtheorem{lemma}[theorem]{Lemma}
\newtheorem{proposition}[theorem]{Proposition}
\newtheorem{corollary}[theorem]{Corollary}
\newtheorem{conjecture}[theorem]{Conjecture}
\theoremstyle{definition}

\newtheorem{example}[theorem]{Example}
\theoremstyle{remark}
\newtheorem{remark}[theorem]{Remark}

\newcommand{\D}{\mathcal D}
\newcommand{\R}{\mathbb R}
\newcommand{\Z}{\mathbb Z}
\newcommand{\avg}[1]{\langle #1\rangle}
\newcommand{\1}{\mathbf 1}
\newcommand{\bmo}{\mathrm{BMO}_d}
\newcommand{\cmo}{\mathrm{CMO}_d}
\newcommand{\eps}{\varepsilon}
\newcommand{\G}{\mathcal G}
\newcommand{\Aw}{[w]_{A_2^d}}
\DeclareMathOperator{\dist}{dist}

\title[$A_2$ growth of dyadic paraproducts]{Scale profile, nesting, and the $A_2$ growth of dyadic paraproducts}
\author{Daewon Chung}
\address{Department of Mathematics, Keimyung University, Daegu, Republic of Korea}
\email{dwchung@kmu.ac.kr}
\date{}

\begin{document}

\begin{abstract}
For the dyadic paraproduct $\pi_b$ and $Q\ge1$, let $\Phi_b(Q)$ be the supremum of $\|\pi_b\|_{L^2(w)}$ over dyadic weights with $[w]_{A_2^d}\le Q$.
Beznosova's linear bound $\Phi_b(Q)\lesssim Q\,\|b\|_{\mathrm{BMO}_d}$ is sharp over $\mathrm{BMO}_d$, while every nonzero symbol safisfies $\Phi_b(Q)\gtrsim Q^{1/2}$. We ask which properties of a fixed symbol place its growth between these extremes, and identify two: the distribution of the Haar coefficients of $b$ over scales and the nesting of their supports. 
Symbols supported on $L$ consecutive scales obey a two-weight bound of order $\sqrt L$; as a consequence, polynomial decay of the scale tails determines the sharp trowth exponent, except at a critical rate where the growth is determined up to a factor $(\log Q)^{1/2}$. For symbols supported on a chain, $\sup_Q Q^{-1/2}\Phi_b(Q)$ is comarable to $\|S_d b\|_{\infty}$, and for every symbol the boundedness of the dyadic square function $S_d b$ is necessary for minimal growth $Q^{1/2}$. Conversely, $\Phi_b(Q)\lesssim Q^{3/4}\|S_d b\|_{\infty}$ for every symbol, and the bound $\|\pi_b\|_{L^2(w)}\lesssim [w]_{A_2^d}^{1/2}\|S_d b\|_{i\infty}$ holds for every symbol at weights whose local $A_2$ imbalance has bounded evergy off a family of dyadic chains with disjoint tops; this class contains all weights used in our lower bounds.
\end{abstract}

\maketitle

\section{Introduction}

\subsection{Background}
Muckenhoupt \cite{Mu72} gave a characterization of the weights $w$ that make the Hardy--Littlewood maximal function bounded on $L^p(w)$. For the Hilbert transform, Hunt, Muckenhoupt and Wheeden \cite{HMW73} established that boundedness is characterized by that same $A_p$ condition.
The quantitative dependence of operator norm on $[w]_{A_p}$ was studied later. The sharp power for the maximal function was obtained by Buckley \cite{Bu93}; its singular integral counterpart came to be called the $A_2$ conjecture. The sharp extrapolation theorem \cite{DGPP05} reduces this problem to establishing linear dependence on $[w]_{A_2}$ in $L^2(w)$.
Dyadic model operators played a decisice role. Wittwer \cite{Wi00} established the linear bound for the martingale transform, while Hukovic, Treil and Volberg \cite{HTV00} obtained it for the dyadic square function. Petermichl \cite{Pe07} proved it for the Hilbert transform; Lacey, Petermichl and Reguera \cite{LPR10}established it for Haar shifts, with bounds polynomial in the complexity obtained in \cite{HPTV}; and Hyt\"onen \cite{Hy12} established it for all Calder\'on--Zygmund operators. Sparse domination subsequently yielded simpler proofs by Lerner \cite{Le13} and Lacey \cite{La17}.
These results were followed by mixed $A_2$--$A_\infty$ refinements \cite{HP13,HL12}.
Let $\D$ denote the standard dyadic intervals of $\R$, $h_I=|I|^{-1/2}(\1_{I_l}-\1_{I_r})$, with $I_l$ and $I_r$ denoting, respectively, the left and right children of $I$, and
\begin{equation}\label{eq:para}
\pi_bf=\sum_{I\in\D}b_I\avg f_Ih_I,\qquad b_I=\langle b,h_I\rangle .
\end{equation}
The dyadic paraproduct forms a basic building block in dyadic harmonic analysis. Indeed, the product $bf$ decomposes into $\pi_bf$, $\pi_b^*f$ and $\pi_fb$; the non-cancellative part of both $T(1)$-type decompositions and commutators \cite{Pe01} is carried by paraproducts. Its $L^2$ norm and $\|b\|_{\bmo}$ are comparable \cite[Theorem~1.32]{Pe01}.
Using Bellman function methods together with a weighted Carleson lemma, Beznosova \cite{Be08} established the linear bound
\begin{equation}\label{eq:bez}
\|\pi_b\|_{L^2(w)}\le C_B\Aw\|b\|_{\bmo},
\end{equation}
and over $\bmo$ this dependence is optimal (Corollary~\ref{cor:window-sharp}).
For commutators $[b,T]$ involving a Calder\'on--Zygmund operator $T$, quadratic dependence is the corresponding sharp dependence \cite{Ch11,CPP12}. Throughout these results, only the $\bmo$ norm of the symbol enters.

\subsection{Symbols with vanishing oscillation}
Coifman, Rochberg and Weiss \cite{CRW76} established that boundedness of commutators with the Riesz transforms is equivalent to $b\in\mathrm{BMO}$. Uchiyama \cite{Uc78} proved that compactness is equivalent to $b$ belonging to $\mathrm{CMO}$, namely the closure of $C_c^\infty$ in $\mathrm{BMO}$, a space refining Sarason's space $\mathrm{VMO}$ \cite{Sa75}. For the dyadic setting, compactness of $\pi_b$ on $L^2$ is equivalent to the vanishing of the mean oscillations of $b$ along small, large and distant dyadic intervals \cite[Theorem~6.2]{HO24}; see also \cite{CP96}. When $p=2$, this condition is equivalent to $b\in\cmo$, defined as the $\bmo$-closure of finite Haar sums.
The essential norm estimate \eqref{eq:ess}, for example, yields compactness on every $L^2(w)$ with $w\in A_2^d$. More generally, compactness extrapolates to weighted spaces \cite{HL23,SVW}.

This leads to a quantitative question: compactness expresses smallness of the coefficients of $b$ at extreme scales and locations, but does that smallness also reduce the growth of $\|\pi_b\|_{L^2(w)}$ as $\Aw\to\infty$?
Two extreme cases point toward an affirmative answer. When the Haar coefficients of $b$ are supported on a single scale, the norm's supremum over $\Aw\le Q$ is comparable to $Q^{1/2}$ (Lemma~\ref{lem:single} and Proposition~\ref{prop:lower}). For the worst symbols in $\bmo$, by contrast, its growth is like $Q$. Because finite Haar sums have $\cmo$ as their closure, intermediate growth might be expected, with the rate of smallness of $b$ at extreme scales determining that growth. We express this precisely by writing $\sigma=w^{-1}$, $\Aw=\sup_{I\in\D}\avg w_I\avg\sigma_I$, and

\begin{equation}\label{eq:Phi}
\Phi_b(Q)=\sup\{\|\pi_b\|_{L^2(w)\to L^2(w)}:\Aw\le Q\},\qquad[b]_{\G_\alpha}=\sup_{Q\ge1}Q^{-\alpha}\Phi_b(Q),
\end{equation}
with $\G_\alpha=\{b:[b]_{\G_\alpha}<\infty\}$. Here the role of $Q\ge1$ is that of a threshold parameter rather than the characteristic of any particular weight; assertions concerning a fixed weight use $\Aw$. The two formulations are equivalent for power bounds: if $\alpha\ge0$ and $K\ge0$ depends only on $b$, then
\begin{equation}\label{eq:QvsW}
\Phi_b(Q)\le C\,Q^\alpha K\ \text{ for all }Q\ge1\quad\Longleftrightarrow\quad\|\pi_b\|_{L^2(w)}\le C\Aw^\alpha K\ \text{ for all }w\in A_2^d,
\end{equation}
since from $\Aw\le Q$ one obtains $\Aw^\alpha\le Q^\alpha$, whereas in the reverse direction one can choose $Q=\Aw\ge1$. Thus $[b]_{\G_\alpha}$ is the smallest constant $C$ in \eqref{eq:QvsW} satisfying $K=1$. In proofs involving a fixed weight, we occasionally use the abbreviation $Q=\Aw$. By testing $\pi_b$ on $\1_J$ with $w=1$, we obtain $\|b\|_{\bmo}\le\Phi_b(1)$; consequently, $\G_1=\bmo$ with $\|b\|_{\bmo}\le[b]_{\G_1}\le C_B\|b\|_{\bmo}$. In the other direction, for each $b$ we have $\Phi_b(Q)\ge(Q/2)^{1/2}\sup_I|b_I||I|^{-1/2}$ (Proposition~\ref{prop:lower}). Hence the minimal growth is $Q^{1/2}$, for $\alpha<\frac12$, only symbols whose Haar coefficients vanish belong to $\G_\alpha$. We ask what properties of $b$ ensure its membership in $\G_\alpha$ for $\frac12\le\alpha<1$.

Write
\begin{equation}\label{eq:tail}
\Delta_kb=\sum_{|I|=2^{-k}}b_Ih_I,\qquad b^{>K}=\sum_{|k|>K}\Delta_kb
\end{equation}
for the component of $b$ at scale $2^{-k}$ and for its scale tail. To describe the scale profile of $b$, we use the decay of $\|b^{>K}\|_{\bmo}$ as $K\to\infty$; to describe its nesting, we use the dyadic square function
\begin{equation}\label{eq:Sd}
S_db=\Big(\sum_{I\in\D}\frac{b_I^2}{|I|}\1_I\Big)^{1/2}.
\end{equation}
For a pair of weights, our joint characteristic is $[u,v]_{A_2^d}=\sup_{I\in\D}\avg{u^{-1}}_I\avg v_I$; consequently, $[w,w]_{A_2^d}=\Aw$.

\subsection{Main results}
Our first result addresses how the distribution varies over scales. A sum of the one-scale estimate across $L$ scales yields a constant of order $L$. For paraproducts, however, the correct order is $\sqrt L$, and the constant is absolute even in the two-weight case. Summation across dyadic blocks of scales then yields sharp exponents away from a critical decay rate.

\begin{theorem}\label{thmA}
\textup{(i)} Suppose the Haar coefficients of $b$ are zero except on $L$ consecutive scales. Then
\[
\|\pi_b\|_{L^2(u)\to L^2(v)}\le C\sqrt{L\,[u,v]_{A_2^d}}\;\|b\|_{\bmo}
\]
where $C$ is an absolute constant. In the one-weight case, taking the supremum of $\|\pi_b\|_{L^2(w)}/\|b\|_{\bmo}$ over $\Aw\le Q$ and all such $b$ gives a quantity comparable to $\sqrt{\min(L,Q)\,Q}$.
\textup{(ii)} Suppose that, for every $K\ge0$, we have $\|b^{>K}\|_{\bmo}\le C_0(1+K)^{-\kappa}$. Then, for $Q\ge2$
\begin{equation}\label{eq:rates}
\Phi_b(Q)\le C_\kappa(\|b\|_{\bmo}+C_0)\cdot\begin{cases}Q^{1-\kappa},&0<\kappa<\frac12,\\Q^{1/2}\log(e+Q),&\kappa=\frac12,\\Q^{1/2},&\kappa>\frac12.\end{cases}
\end{equation}
Given any $0<\kappa<\frac12$, one can find a bounded $b\in\cmo$ satisfying $\|b^{>K}\|_{\bmo}\le\sqrt2(1+K)^{-\kappa}$ for every $K$ as well as $\Phi_b(Q)\ge\frac1{\sqrt8}(Q/2)^{1-\kappa}$ for every $Q\ge2$.
\end{theorem}

The logarithm is necessary at the critical rate $\kappa=\frac12$: a bounded symbol belonging to $\cmo$, having this tail decay and satisfying $\Phi_b(Q)\simeq(Q\log Q)^{1/2}$, is provided by Corollary~\ref{cor:critical}.
It remains open whether $\frac12$ can replace the power of the logarithm in Theorem~\ref{thmA}(ii). An affirmative answer to question \eqref{eq:Sq} below would imply this improvement (Remark~\ref{rem:critical}).

The quantity $\Phi_b$ cannot be determined from the scale profile alone: symbols having one coefficient of equal size on each of $L$ scales may exhibit growth $Q^{1/2}$ when their supports are disjoint, or growth at least $L/\sqrt8$ at $Q=L$ when those supports are nested (Proposition~\ref{prop:notdelta}). For fully nested supports forming a chain, weighted Hardy operators give the reduction of the operator; the appropriate quantity is then an $\ell^2$ sum over that chain, namely $\|S_db\|_\infty$.
The square function is, in general, necessary for minimal growth.
\begin{theorem}\label{thmB}
\textup{(i)} Let $I_0\supsetneq I_1\supsetneq\cdots$, where $I_{k+1}$ is a child of $I_k$, and let $b=\sum_kc_k|I_k|^{1/2}h_{I_k}$. Then
\begin{equation}\label{eq:chainB}
\|\pi_b\|_{L^2(u)\to L^2(v)}\le8\,[u,v]_{A_2^d}^{1/2}\Big(\sum_kc_k^2\Big)^{1/2},\qquad[b]_{\G_{1/2}}\ge\frac1{\sqrt{12}}\Big(\sum_kc_k^2\Big)^{1/2},
\end{equation}
and $(\sum_kc_k^2)^{1/2}=\|S_db\|_\infty$.
\textup{(ii)} With an absolute constant $C$, one has $\|S_db\|_\infty\le C[b]_{\G_{1/2}}$ for every $b$. Consequently, $\G_{1/2}\subseteq\{b:S_db\in L^\infty\}$.
\end{theorem}

The proofs yield the constants $8$ and $1/\sqrt{12}$ appearing in \eqref{eq:chainB}; their optimality is not asserted. In the other statements, upper bounds use absolute constants $C$ unless a named constant is specified, such as $C_B$ in \eqref{eq:bez}. Likewise, the explicit numbers occurring in proofs, technical lemmas, and lower-bound constructions arise from the arguments, with no assertion of optimality. The signs of the coefficients have no effect on the lower bounds of Theorem~\ref{thmB}. These bounds are obtained by testing the adjoint against a weight singular at the bottom of a chain: the contribution of that chain is a causal Toeplitz operator with a symbol that does not vanish on the unit circle. When a chain has nonincreasing $|c_k|\to0$, combining the upper and lower bounds gives the exact growth $\Phi_b(Q)\simeq Q^{1/2}(\sum_{k<Q}c_k^2)^{1/2}$ (Theorem~\ref{thm:monotone}).

The role of compactness is also clarified by these results. Although $b\in\cmo$ implies $\Phi_b(Q)=o(Q)$ (Corollary~\ref{cor:qual}), no universal rate below linear growth follows from compactness alone, even for bounded symbols belonging to the classical space $\mathrm{CMO}(\R)$ (Theorem~\ref{thm:norate}). By contrast, minimal growth $\Phi_b(Q)=Q^{1/2}$ is possible for a symbol whose distance from $\cmo$ is one; only in the worst case does the distance to $\cmo$ control the essential norm (Theorem~\ref{thm:ess}).

\subsection{Symbols with bounded square function}\label{sec:Sq}
By Theorem~\ref{thmB}, the minimal growth $Q^{1/2}$ requires boundedness of $S_db$; on chains, this condition is sufficient as well. We establish two unconditional results in the converse direction. For arbitrary weights and every symbol,
\begin{equation}\label{eq:34intro}
\|\pi_b\|_{L^2(w)}\le C\Aw^{3/4}\|S_db\|_\infty
\end{equation}
(Proposition~\ref{prop:34}), which yields $\{b:S_db\in L^\infty\}\subseteq\G_{3/4}$. The exponent $\frac34$ admits replacement by $\frac12$ for weights whose local $A_2^d$ imbalance has bounded energy off a chain of dyadic intervals, or off a family of chains with disjoint tops (Theorems~\ref{thm:spine} and~\ref{thm:forest}). Included in this class are the power weights and all weights appearing in this paper’s lower bounds (Corollary~\ref{cor:power}, Remark~\ref{rem:blocks}). It remains open whether
\begin{equation}\label{eq:Sq}
\|\pi_b\|_{L^2(w)}\le C\Aw^{1/2}\|S_db\|_\infty\qquad\text{for all $b$ and all $w\in A_2^d$},
\end{equation}
equivalently, whether $\G_{1/2}=\{b:S_db\in L^\infty\}$ with $[b]_{\G_{1/2}}\simeq\|S_db\|_\infty$. There are three observations that support \eqref{eq:Sq} as the natural candidate for a weight-free description of $\G_{1/2}$.

First, nested Haar supports generate all lower bounds beyond $Q^{1/2}$ obtained in this paper. Chains underlie the examples for Corollary~\ref{cor:window-sharp}, Proposition~\ref{prop:rate-sharp}, Corollary~\ref{cor:critical} and Theorem~\ref{thm:norate}; by contrast, Proposition~\ref{prop:notdelta} establishes minimal growth when those same coefficients are put on disjoint intervals. On a chain, Theorem~\ref{thmB}(i) identifies $[b]_{\G_{1/2}}$, within absolute constants, as the $\ell^2$ sum of the normalized coefficients along that chain. For a general symbol, $S_db(x)^2$ likewise denotes the $\ell^2$ sum along the dyadic path that terminates at $x$. Consequently, \eqref{eq:Sq} would assert that no branching Haar support is worse than its worst path.

Second, the upper bounds currently available agree with \eqref{eq:Sq}. By Lemma~\ref{lem:single}, the norms of the one-scale pieces $\pi_{\Delta_kb}$ do not exceed $Q^{1/2}\|\Delta_kb\|_\infty$. Their sum yields the $\ell^1$ bound in Corollary~\ref{cor:qual}(i), while square summation along paths corresponds to $\|S_db\|_\infty\le(\sum_k\|\Delta_kb\|_\infty^2)^{1/2}$. Likewise, support of $b$ on $L$ consecutive scales implies $\|S_db\|_\infty\le\sqrt L\,\|b\|_{\bmo}$, according to Lemma~\ref{lem:basic}(i). Thus the one-weight case of the bound $C\sqrt{L\Aw}\,\|b\|_{\bmo}$ from Theorem~\ref{thmA}(i) would be contained in \eqref{eq:Sq}.

Third, the loss arising in the weighted Carleson embedding is eliminated by the hypothesis $S_db\in L^\infty$. Lemma~\ref{lem:univ} below identifies $\|S_db\|_\infty^2$ as the best constant, uniformly over all measures, in the Carleson condition for the coefficients of $b$. In contrast, the analogous condition for $b\in\bmo$ with respect to $\sigma\,dx$ requires an additional factor $[\sigma]_{A_\infty}$ (Remark~\ref{rem:partial}). Thus, for every weight, the component of $\pi_b$ orthogonal in $L^2(w)$ obeys \eqref{eq:Sq} (Proposition~\ref{prop:drift}). The remaining component is dominated by a positive dyadic operator with coefficients given by $|b_I||I|^{-1/2}$ times the relative imbalance of $w$ between the two children of $I$. Whether this imbalance incurs another power of $\Aw$ is the question.

We state the anticipated answer in Section~\ref{sec:concl} as Conjecture~\ref{conj:S} and explain there at which point the proof of \eqref{eq:34intro} loses the factor $\Aw^{1/4}$.

\subsection*{Organization} Section~\ref{sec:prelim} introduces the notation along with two tools: a discrete Hardy inequality and a two-weight testing theorem for positive dyadic operators. The universal lower bound is proved in Section~\ref{sec:single}, which addresses one scale. The proof of Theorem~\ref{thmA} occupies Sections~\ref{sec:window} and~\ref{sec:rates}. Section~\ref{sec:tree} proves Theroem~\ref{thmB} and contains the result on nesting and chains and on symbols with bounded square function. Compactness is addressed in Sections~\ref{sec:norate} and~\ref{sec:compact}; the conjecture and concluding remarks appear in Section~\ref{sec:concl}.

\section{Preliminaries}\label{sec:prelim}

Throughout, we take $\pi_b$, $\Phi_b$, $b^{>K}$ and $S_db$ as specified in \eqref{eq:para}, \eqref{eq:Phi}, \eqref{eq:tail} and \eqref{eq:Sd}. For $I\in\D$, put $\D_k=\{I:|I|=2^{-k}\}$, $\avg f_I=|I|^{-1}\int_If$ and $f_I=\langle f,h_I\rangle$. When $x\in I$, we write $f_Ih_I(x)=\avg f_{I'}-\avg f_I$, with $I'$ denoting the child of $I$ that contains $x$. We define
\begin{equation}\label{eq:bmo}
\|\Delta_kb\|_\infty=\sup_{I\in\D_k}\frac{|b_I|}{|I|^{1/2}},\qquad\|b\|_{\bmo}^2=\sup_{J\in\D}\frac1{|J|}\sum_{I\subseteq J}b_I^2 ,
\end{equation}
and we use the notation $b_{\mathcal F}=\sum_{I\in\mathcal F}b_Ih_I$ and $b^S=\sum_{k\in S}\Delta_kb$ for $\mathcal F\subseteq\D$ and $S\subseteq\Z$. Symbols sharing the same Haar coefficients are identified, because $\pi_b$ depends solely on the coefficients $(b_I)$; on functions, $\|b\|_{\bmo}$, $\Phi_b$ and $[b]_{\G_\alpha}$ are seminorms. The $\bmo$ norm cannot increase when the coefficients are restricted to a subfamily. From $\frac1{|J|}\sum_{I\subseteq J}b_I^2=\frac1{|J|}\int_J\sum_{I\subseteq J}\frac{b_I^2}{|I|}\1_I$ it follows that $\|b\|_{\bmo}\le\|S_db\|_\infty$.
\begin{lemma}\label{lem:basic}
\textup{(i)} $\|\Delta_kb\|_\infty\le\|b\|_{\bmo}$, with $\|\Delta_kb\|_\infty\le\|b^{>K}\|_{\bmo}$ for $|k|>K$. \textup{(ii)} $\|b^S\|_{\bmo}^2\le\sum_{k\in S}\|\Delta_kb\|_\infty^2$. \textup{(iii)} Whenever $b\in\cmo$, we obtain $\|b^{>K}\|_{\bmo}\to0$.
\end{lemma}

\begin{proof}
(i) From \eqref{eq:bmo} we obtain $b_I^2\le\sum_{I'\subseteq I}b_{I'}^2\le\|b\|^2_{\bmo}|I|$; use this for both $b$ and its scale tail. (ii) Once $J$ and $k$ are fixed, those intervals of $\D_k$ that lie in $J$ are disjoint. Their contribution to the Carleson sum over $J$ therefore does not exceed $\|\Delta_kb\|_\infty^2|J|$. (iii) Suppose $\|b-b_{\mathcal F}\|_{\bmo}<\eta$ and $\mathcal F$ is finite. For sufficiently large $K$, the scale tails of $b$ coincide with those of $b-b_{\mathcal F}$, giving $\|b^{>K}\|_{\bmo}\le\|b-b_{\mathcal F}\|_{\bmo}<\eta$.
\end{proof}
A weight $w$ is positive a.e. and satisfies $w,\sigma=w^{-1}\in L^1_{\mathrm{loc}}$. The dyadic Fujii--Wilson constant is given by
\begin{equation}\label{eq:Ainf}
[w]_{A_\infty}=\sup_{J\in\D}\frac1{w(J)}\int_JM^d(w\1_J)\le e\,\Aw
\end{equation}
as follows from the dyadic version of \cite[Proposition~2.2]{HP13}. Given two weights, set $\sigma=u^{-1}$; the meaning of $[u,v]_{A_2^d}<\infty$ is then $\sigma(I)v(I)\le [u,v]_{A_2^d}|I|^2$ for every $I\in\D$. Absolute convergence of the series $\langle\pi_bf,g\rangle=\sum_Ib_I\avg f_Ig_I$ for $b\in\bmo$ and $f,g\in L^\infty_c$ allows norms to be computed on these dense classes; moreover, the pairing is additive when the coefficients are split. In particular, suppose an operator bound holds for each finite subfamily of the coefficients of $b$, with a constant that cannot increase upon restriction (for example, $\|b\|_{\bmo}$ or $\sum_kc_k^2$ below). Dominated convergence in the pairing then yields the same bound for $b$. Passages from finite coefficient families to infinite ones will use this observation without further mention. A family $\mathcal S\subseteq\D$ is sparse when pairwise disjoint sets $E(S)\subseteq S$ exist satisfying $|E(S)|\ge|S|/2$; consequently, $\sum_{S'\in\mathcal S,S'\subseteq S}|S'|\le2|S|$.
\begin{theorem}[{\cite[Theorems~1.6 and 1.11]{LSU}, \cite[Theorem~1.2]{HHL}}]\label{thm:LSU}
Suppose that $\lambda_I\ge0$ and $T(f\sigma)=\sum_{I\in\D}\lambda_I\avg{f\sigma}_I\1_I$. We then have $\|T(\cdot\,\sigma)\|_{L^2(\sigma)\to L^2(v)}\le C_T(\mathfrak T+\mathfrak T^*)$ for an absolute constant $C_T$, where
\begin{equation}\label{eq:testing}
\mathfrak T^2=\sup_{R\in\D}\frac1{\sigma(R)}\int\Big(\sum_{I\subseteq R}\lambda_I\avg\sigma_I\1_I\Big)^2dv
\end{equation}
and the definition of $\mathfrak T^*$ is obtained by interchanging the roles of $\sigma$ and $v$.
\end{theorem}

In dimension one, the characterization for $p=q=2$ is due to Nazarov, Treil and Volberg \cite{NTV99}. The normalization of the operator in \cite{HHL} uses $\int_I f\,d\sigma$ rather than $\avg{f\sigma}_I$; equivalently, $\lambda_I$ is replaced by $\lambda_I/|I|$.

\begin{lemma}\label{lem:hardy}
Consider $u_j,v_j,x_j\ge0$, with indices ranging over a set of consecutive integers, adopting the convention $x_i^2/v_i=0$ when $x_i=v_i=0$ and $x_i^2/v_i=\infty$ when $x_i>0=v_i$. It follows that
\begin{equation}\label{eq:hardy}
\sum_ju_j\Big(\sum_{i\ge j}x_i\Big)^2\le4B_1\sum_i\frac{x_i^2}{v_i},\qquad\sum_ju_j\Big(\sum_{i<j}x_i\Big)^2\le4B_2\sum_i\frac{x_i^2}{v_i},
\end{equation}
where we set $B_1=\sup_n(\sum_{j\le n}u_j)(\sum_{i\ge n}v_i)$ and $B_2=\sup_n(\sum_{j>n}u_j)(\sum_{i\le n}v_i)$.
\end{lemma}

\begin{proof}
When the index set is infinite, use finite sections and then monotone convergence: the constants $B_1,B_2$ for these sections are no larger than those for the full sequences. When the index set is finite, substitute $u_j+\epsilon,v_j+\epsilon$ for $u_j,v_j$ and take $\epsilon\to0$. Both the constants and the right-hand sides converge, the latter according to the stated convention. We can therefore restrict attention to finitely many indices with $u_j,v_j>0$. Set $V_j=\sum_{i\ge j}v_i$ and $U_j=\sum_{l\le j}u_l$. Comparison with $\int t^{-1/2}dt$, using the decrease of $V_i$ and the increase of $U_j$, yields $\sum_{i\ge j}v_iV_i^{-1/2}\le2V_j^{1/2}$ and $\sum_{j\le i}u_jU_j^{-1/2}\le2U_i^{1/2}$. The Cauchy--Schwarz inequality now gives
\[
\Big(\sum_{i\ge j}x_i\Big)^2\le\sum_{i\ge j}\frac{x_i^2V_i^{1/2}}{v_i}\sum_{i\ge j}v_iV_i^{-1/2}\le2V_j^{1/2}\sum_{i\ge j}\frac{x_i^2V_i^{1/2}}{v_i}.
\]
The left-hand side of the first inequality is bounded above by $2B_1^{1/2}\sum_i(x_i^2V_i^{1/2}/v_i)\sum_{j\le i}u_jU_j^{-1/2}\le4B_1^{1/2}\sum_i(x_i^2/v_i)(U_iV_i)^{1/2}\le4B_1\sum_ix_i^2/v_i$ after applying $V_j\le B_1/U_j$ and interchanging the sums. Reversing the order of the indices gives the second inequality by the same argument.
\end{proof}

\section{One scale}\label{sec:single}

\begin{lemma}\label{lem:single}
Suppose $\mathcal A\subseteq\D$ is a collection of pairwise disjoint intervals and $Tf=\sum_{I\in\mathcal A}c_I\avg f_Ih_I$. In that case, $T:L^2(u)\to L^2(v)$ is an orthogonal sum of rank-one operators whose norms are
\begin{equation}\label{eq:sI}
s_I=|c_I||I|^{-1/2}(\avg\sigma_I\avg v_I)^{1/2}.
\end{equation}
Consequently, $\|T\|=\sup_Is_I$. Moreover, the nonzero singular values of $T$ are precisely the numbers $s_I$, and $\|T\|_{\mathrm{ess}}=\inf_{\mathcal F}\sup_{I\in\mathcal A\setminus\mathcal F}s_I$, where the infimum is over finite $\mathcal F$. In particular,
\[
\|\pi_{\Delta_kb}\|_{L^2(u)\to L^2(v)}\le\|\Delta_kb\|_\infty[u,v]_{A_2^d}^{1/2}.
\]
\end{lemma}

\begin{proof}
The values of the piece $f\mapsto c_I\avg f_Ih_I$ are functions supported on $I$, while its dependence is solely on $f\1_I$; hence the pieces act between mutually orthogonal subspaces. The Cauchy--Schwarz inequality gives its norm as $|c_I|\,\|h_I\|_{L^2(v)}\sup\{|\avg f_I|:\|f\1_I\|_{L^2(u)}\le1\}=|c_I|\avg v_I^{1/2}|I|^{-1}\sigma(I)^{1/2}$, with equality attained at $f=\sigma\1_I$.
\end{proof}

\begin{proposition}\label{prop:lower}
For every $b$ and $Q\ge1$, $\Phi_b(Q)\ge(Q/2)^{1/2}\sup_k\|\Delta_kb\|_\infty$.
\end{proposition}

\begin{proof}
Fix $I_0\in\D$ and set $a=2Q-1\ge1$; also let $I_{0+}=(I_0)_l$, $I_{0-}=(I_0)_r$ and $w_\pm=a\1_{I_{0\pm}}+\1_{\R\setminus I_{0\pm}}$.
 A dyadic interval $J$ has three possibilities: it avoids the interior of $I_{0\pm}$, lies inside that interior, yielding $\avg{w_\pm}_J\avg{\sigma_\pm}_J=1$, or contains $I_0$, yielding $\avg{w_\pm}_J\le\frac{1+a}2=Q$ and $\avg{\sigma_\pm}_J\le1$.
Consequently, $[w_\pm]_{A_2^d}\le Q$ as well as $\avg{w_\pm}_{I_0}\avg{\sigma_\pm}_{I_0}=\frac{(a+1)^2}{4a}\ge Q/2$. Use $f=\sigma_\pm\1_{I_0}$ and $g=w_\pm h_{I_0}$ as tests.
For $I\subsetneq I_0$, constancy of the function $g$ on $I$ gives $g_I=0$. The contribution from the interval $I_0$ is $b_{I_0}\avg{\sigma_\pm}_{I_0}\avg{w_\pm}_{I_0}$, identical for either sign, whereas intervals disjoint from $I_0$ give no contribution. The ancestors $I\supsetneq I_0$ supply $X\int_{I_0}h_{I_0}w_\pm$; here $X=\sum_{I\supsetneq I_0}b_I\avg{f}_Ih_I(I_0)$ is independent of the sign since $\sigma_\pm(I_0)$ is independent of it, and $\int_{I_0}h_{I_0}w_\pm=\pm|I_0|^{1/2}\frac{a-1}2$. Taking the sign that prevents cancellation between the two contributions and applying $\|f\|_{L^2(w_\pm)}=\sigma_\pm(I_0)^{1/2}$ together with $\|g\|_{L^2(\sigma_\pm)}=\avg{w_\pm}_{I_0}^{1/2}$ yields $\|\pi_b\|_{L^2(w_\pm)}\ge|b_{I_0}||I_0|^{-1/2}(\avg{w_\pm}_{I_0}\avg{\sigma_\pm}_{I_0})^{1/2}\ge(Q/2)^{1/2}|b_{I_0}||I_0|^{-1/2}$.
\end{proof}

\begin{corollary}\label{cor:trivial}
Under the conditions $\alpha<\frac12$ and $b\in\G_\alpha$, every Haar coefficient of $b$ vanishes.
\end{corollary}

\section{Windows of scales}\label{sec:window}

Let $W$ comprise $L$ consecutive integers, $\D_W=\bigcup_{k\in W}\D_k$; say that $b$ is $W$-supported when $b_I=0$ for $I\notin\D_W$.
\begin{lemma}\label{lem:sparse}
Whenever $b$ is $W$-supported and $f,g\in L^\infty_c$, a sparse family $\mathcal S\subseteq\D_W$ exists with
\begin{equation}\label{eq:sparse}
|\langle\pi_bf,g\rangle|\le36\|b\|_{\bmo}\sum_{S\in\mathcal S}\avg{|f|}_S\avg{|g|}_S|S| .
\end{equation}
\end{lemma}

\begin{proof}
The coarsest scale $k_0$ intervals of $W$ give a partition of $\R$ in which each $I\in\D_W$ belongs to precisely one interval. Since just finitely many such intervals intersect the supports of $f$ and $g$, the pairing decomposes accordingly; fix any one of these intervals $J_0$. Beginning with $S=J_0$, define the stopping children of $S$ to be the maximal $I\in\D_W$ satisfying $I\subsetneq S$ and either $\avg{|f|}_I>4\avg{|f|}_S$ or $\avg{|g|}_I>4\avg{|g|}_S$. Their total measure does not exceed $|S|/2$; consequently, $E(S)=S\setminus\bigcup\{\text{stopping children}\}$ has $|E(S)|\ge|S|/2$. Iteration of this construction produces a sparse family $\mathcal S$ of stopping intervals. Define $\mathcal T(S)$ as the collection of intervals in $\D_W$ having $S$ as their minimal ancestor in $\mathcal S$. The averages of $|f|$ and $|g|$ on $\mathcal T(S)$ are bounded by $4$ times their averages on $S$; on children of intervals in $\mathcal T(S)$, the corresponding bound is $8$ times their averages on $S$. For $x\in S$, consecutiveness of the scales in $W$ implies that the intervals in $\mathcal T(S)$ containing $x$ constitute a chain of consecutive generations beginning at $S$. Its final interval has a child containing $x$ that is either a stopping interval or outside $\D_W$. Along this chain, telescoping yields $G_S(x)=\avg g_{I'_x}-\avg g_S$ for $G_S=\sum_{I\in\mathcal T(S)}g_Ih_I$, with $I'_x$ denoting the child containing $x$ of the smallest interval in $\mathcal T(S)$ that contains $x$. Thus $|G_S|\le9\avg{|g|}_S\1_S$ and $\sum_{\mathcal T(S)}g_I^2=\|G_S\|_2^2\le81\avg{|g|}_S^2|S|$. It follows that
\[
\Big|\sum_{I\in\mathcal T(S)}b_I\avg f_Ig_I\Big|\le4\avg{|f|}_S\Big(\sum_{I\subseteq S}b_I^2\Big)^{1/2}\Big(\sum_{I\in\mathcal T(S)}g_I^2\Big)^{1/2}\le36\|b\|_{\bmo}\avg{|f|}_S\avg{|g|}_S|S|,
\]
and the lemma follows upon summing over $S\in\mathcal S$.
\end{proof}

\begin{theorem}\label{thm:window}
Suppose $b$ is $W$-supported and satisfies $|W|=L$. We then have
\begin{equation}\label{eq:window2}
\|\pi_b\|_{L^2(u)\to L^2(v)}\le C\sqrt{L\,[u,v]_{A_2^d}}\;\|b\|_{\bmo},
\end{equation}
and in the one weight case
\begin{equation}\label{eq:window1}
\|\pi_b\|_{L^2(w)}\le C\Aw^{1/2}\Big(\min(L,2[\sigma]_{A_\infty})^{1/2}+\min(L,2[w]_{A_\infty})^{1/2}\Big)\|b\|_{\bmo},
\end{equation}
with an absolute constant $C$.
\end{theorem}
\begin{proof}
Set $f=\varphi\sigma$ and $g=\psi v$; hence $\|f\|_{L^2(u)}=\|\varphi\|_{L^2(\sigma)}$ and $\|g\|_{L^2(v^{-1})}=\|\psi\|_{L^2(v)}$. We can express the sparse form in \eqref{eq:sparse} as $\langle A_{\mathcal S}(|\varphi|\sigma),|\psi|v\rangle$, with $A_{\mathcal S}h=\sum_{S\in\mathcal S}\avg h_S\1_S$ a positive dyadic operator satisfying $\lambda_S=\1_{\mathcal S}(S)$. For $R\in\D$, by sparseness together with $(\sum_Sa_S\1_S)^2\le2\sum_Sa_S\1_S\sum_{S'\subseteq S}a_{S'}\1_{S'}$ and $\avg\sigma_{S'}v(S')\le [u,v]_{A_2^d}|S'|$,

\begin{equation}\label{eq:testcalc}
\int\Big(\sum_{S\subseteq R}\avg\sigma_S\1_S\Big)^2dv\le2\sum_{S\subseteq R}\avg\sigma_S\sum_{S'\subseteq S}\avg\sigma_{S'}v(S')\le4[u,v]_{A_2^d}\sum_{S\subseteq R}\sigma(S).
\end{equation}
Because intervals belonging to the same scale are disjoint and there are at most $L$ scales, $\sum_{S\subseteq R}\sigma(S)\le L\sigma(R)$. Moreover, $\sum_{S\subseteq R}\sigma(S)\le2\sum_{S\subseteq R}\avg\sigma_S|E(S)|\le2\int_RM^d(\sigma\1_R)\le2[\sigma]_{A_\infty}\sigma(R)$. Thus, by \eqref{eq:testcalc}, the testing constant \eqref{eq:testing} of $A_{\mathcal S}$ obeys $\mathfrak T\le2([u,v]_{A_2^d}\min(L,2[\sigma]_{A_\infty}))^{1/2}$; the symmetric conclusion holds for $\mathfrak T^*$. Combining \eqref{eq:sparse} with Theorem~\ref{thm:LSU}, valid uniformly over sparse families, yields \eqref{eq:window2} with $C=144\,C_T$ as well as \eqref{eq:window1} with $C=72\,C_T$.
\end{proof}

\begin{proposition}\label{prop:chain}
Let $P\in\D$ have left-descendant chain $P=I_0\supsetneq I_1\supsetneq\cdots$, and take $0<\delta\le\frac12$, $\eps>0$, $L\ge1$, and

\begin{equation}\label{eq:chainex}
b=\eps\sum_{k<L}(-1)^k|I_k|^{1/2}h_{I_k},\qquad w_\delta=\Big(\frac{x-\inf P}{|P|}\Big)^{1-\delta}\1_P+\1_{\R\setminus P}.
\end{equation}
Then we have $\|b\|_{\bmo}\le\sqrt2\eps$, $\|b\|_\infty\le2\eps$, and $[w_\delta]_{A_2^d}\le1/\delta$, while $f=\sigma\1_P$ obeys
\begin{equation}\label{eq:chainlower}
\|\pi_bf\|_{L^2(w_\delta)}\ge\eps\sqrt{\min(L,1/\delta)/(8\delta)}\;\|f\|_{L^2(w_\delta)} .
\end{equation}
\end{proposition}

\begin{proof}
All quantities remain invariant under the affine change of variables mapping $P$ to $[0,1)$. We therefore present the computation for $P=[0,1)$, $I_k=[0,2^{-k})$; the lengths of the ancestors of $P$ are all that enter below. For $b$, the Carleson sums are $\eps^2\sum_{I_k\subseteq J}|I_k|\le2\eps^2|J|$. On the right child $A_j=I_j\setminus I_{j+1}$ of $I_j$, $b=\eps[\sum_{k<j}(-1)^k-(-1)^j]$ holds, yielding $|b|\le2\eps$. When $J=I_m$, we obtain $\avg w_J\avg\sigma_J=\frac1{\delta(2-\delta)}$. Any dyadic $J\subseteq[0,1)$ that does not contain $0$ obeys $\inf J\ge|J|$. Consequently, the variation of $w$ on $J$ is by a factor at most $2^{1-\delta}\le2$, and $\avg w_J\avg\sigma_J\le2$. When $J\supsetneq[0,1)$, both $\avg w_J\le1$ and $\avg\sigma_J\le\frac12(1+\frac1\delta)$ hold. All of these are at most $1/\delta$ because $\delta\le\frac12$. Taking $r=2^{1-\delta}$ gives $\avg f_{I_k}=\avg\sigma_{I_k}=r^k/\delta$; hence, on $A_j$, $j<L$,

\begin{equation}\label{eq:pointwise}
|\pi_bf|=\frac\eps\delta\Big|\sum_{k<j}(-r)^k-(-r)^j\Big|=\frac\eps\delta\Big|(-r)^j-\frac{1-(-r)^j}{1+r}\Big|\ge\frac\eps\delta\,r^j,
\end{equation}
as $r^j(1+\frac1{1+r})-\frac1{1+r}\ge r^j$. Squaring \eqref{eq:pointwise} and then summing over the integers $0\le j<\min(L,1/\delta)$ yields $\|\pi_bf\|^2\ge\frac{\eps^2}{8\delta^2}\min(L,\frac1\delta)$, since $w(A_j)=2^{-(j+1)(2-\delta)}\frac{2^{2-\delta}-1}{2-\delta}\ge2^{-(j+1)(2-\delta)}$ and, for $j<1/\delta$, $r^{2j}2^{-j(2-\delta)}=2^{-j\delta}\ge\frac12$; meanwhile, $\|f\|^2=\sigma(P)=\frac1\delta$.
\end{proof}

With equal signs, near cancellation occurs between the two terms in \eqref{eq:pointwise} as $r\to2$. Alternation is therefore necessary for this test pair, whereas the signs have no role for the dual pair of Proposition~\ref{prop:dualchain}.

\begin{corollary}\label{cor:window-sharp}
Given $Q\ge2$ and $L\ge1$,
\begin{equation}\label{eq:MLQ}
\tfrac14\sqrt{\min(L,Q)\,Q}\le M_{L,Q}\le C_W\sqrt{\min(L,Q)\,Q},
\end{equation}
where the supremum of $\|\pi_b\|_{L^2(w)}/\|b\|_{\bmo}$ over both $W$-supported $b$ satisfying $|W|=L$ and $\Aw\le Q$ is denoted by $M_{L,Q}$.
In particular, over $\bmo$, Beznosova's bound is sharp.
\end{corollary}

\begin{proof}
For the lower bound, choose $\delta=1/Q$ in \eqref{eq:chainex} and \eqref{eq:chainlower}, then apply $\|b\|_{\bmo}\le\sqrt2\eps$. Combining \eqref{eq:window1} with $[w]_{A_\infty},[\sigma]_{A_\infty}\le c\Aw$ yields the upper bound. The choice $L\ge Q$ yields the lower bound $Q/4$.
\end{proof}

\section{Vanishing rates}\label{sec:rates}

\begin{corollary}\label{cor:qual}
\textup{(i)} $\Phi_b(Q)\le Q^{1/2}\sum_k\|\Delta_kb\|_\infty$. \textup{(ii)} Provided that $\|b^{>K}\|_{\bmo}\to0$, and in particular when $b\in\cmo$, we have $\Phi_b(Q)=o(Q)$.
\end{corollary}

\begin{proof}
(i) Under the assumption $\sum_k\|\Delta_kb\|_\infty<\infty$, Lemma~\ref{lem:single} gives $\pi_b=\sum_k\pi_{\Delta_kb}$, with convergence in operator norm. (ii) Decompose $b$ at scale $K$. The contribution from the scales $|k|\le K$ is at most $(2K+1)\|b\|_{\bmo}Q^{1/2}$ by (i) and Lemma~\ref{lem:basic}(i), whereas \eqref{eq:bez} bounds the tail contribution by $C_BQ\|b^{>K}\|_{\bmo}$. Consequently, for every $K$ we have $\limsup_{Q\to\infty}\Phi_b(Q)/Q\le C_B\|b^{>K}\|_{\bmo}$.
\end{proof}

\begin{theorem}\label{thm:rates}
Whenever $\|b^{>K}\|_{\bmo}\le C_0(1+K)^{-\kappa}$ for every $K\ge0$, we have \eqref{eq:rates}.
\end{theorem}

\begin{proof}
Let $N=\lceil Q\rceil$. Decompose the scales into $\{0\}$, the blocks $B_j^\pm=\pm\{k:2^j\le k<2^{j+1},\ k\le N\}$ indexed by $0\le j\le\log_2N$, and the tail $\{|k|>N\}$. Since every block contains at most $2^j$ consecutive scales lying in $\{|k|\ge2^j\}$, we obtain $\|b^{B_j^\pm}\|_{\bmo}\le\|b^{>2^j-1}\|_{\bmo}\le C_0\,2^{-j\kappa}$. The upper bound in \eqref{eq:MLQ} therefore bounds the paraproduct of each block by $C_WC_0\,Q^{1/2}2^{j(1/2-\kappa)}$.
By \eqref{eq:sI}, the cost of the scale $0$ is at most $Q^{1/2}\|b\|_{\bmo}$; meanwhile, \eqref{eq:bez} bounds the tail by $C_BC_0\,Q(1+N)^{-\kappa}\le C_BC_0\,Q^{1-\kappa}$. The three cases follow upon summing the geometric series in $j\le\log_2N$.
\end{proof}

\begin{proposition}\label{prop:rate-sharp}
Let $0<\kappa<\frac12$, $P_n=[n,n+1)$, and $\eta_n=2^{-n\kappa}$. Define $b$ by summing over $n\ge1$ the alternating chains from Proposition~\ref{prop:chain} on $P_n$ with $\eps=\eta_n$ and $L=2^n$. We then have $b\in\cmo\cap L^\infty$, together with $\|b^{>K}\|_{\bmo}\le\sqrt2(1+K)^{-\kappa}$ for every $K\ge0$ and $\Phi_b(Q)\ge\frac1{\sqrt8}(Q/2)^{1-\kappa}$ for $Q\ge2$.
\end{proposition}

\begin{proof}
The blocks occupy disjoint unit intervals, with their scales restricted to $0\le k<2^n$. For $J\in\D$, any subfamily of coefficients has Carleson sum over $J$ bounded above by $2\eta^2|J|$; here $\eta$ denotes the largest $\eta_n$ of the participating blocks. When $J$ is contained in one $P_n$, the single chain bound yields this conclusion; otherwise, sum over the unit intervals contained in $J$.
Only blocks satisfying $2^n>K+1$ contribute at scales $k>K$. For these blocks $\eta_n<(1+K)^{-\kappa}$, which yields $\|b^{>K}\|_{\bmo}\le\sqrt2(1+K)^{-\kappa}$.
Likewise, the $\bmo$ norm of the blocks satisfying $n\ge N$ is bounded above by $\sqrt2\eta_N\to0$. The finitely many blocks left over are finite Haar sums; consequently $b\in\cmo$, and moreover $|b|\le2\sup_n\eta_n$.
No coefficients occur on the dyadic ancestors of $P_n$, whereas all other blocks are disjoint from $P_n$. Thus, for $f=\sigma\1_{P_n}$, the function $\pi_bf$ agrees with the paraproduct associated with the $n$-th block.
Applying estimate \eqref{eq:chainlower} with $\delta=2^{-n}$ yields $\Phi_b(2^n)\ge\eta_n2^n/\sqrt8=2^{n(1-\kappa)}/\sqrt8$. When $2^n\le Q<2^{n+1}$, the claim follows by monotonicity of $\Phi_b$.
\end{proof}

\begin{proof}[Proof of Theorem~\ref{thmA}]
Part (i) is a consequence of Theorem~\ref{thm:window} together with Corollary~\ref{cor:window-sharp}; part (ii) follows by combining Theorem~\ref{thm:rates} with Proposition~\ref{prop:rate-sharp}.
\end{proof}

\begin{remark}\label{rem:rates}
For $\kappa>\frac12$, Proposition~\ref{prop:lower} shows that every symbol having a nonzero coefficient attains the exponent $\frac12$ of Theorem~\ref{thmA}(ii).
For $\kappa=\frac12$, a logarithmic loss of at least power $\frac12$ is necessary (Corollary~\ref{cor:critical}).
Under the assumptions $\|\Delta_kb\|_\infty\lesssim(1+|k|)^{-\alpha}$ and $\frac12<\alpha<1$, Lemma~\ref{lem:basic}(ii) yields $\|b^{>K}\|_{\bmo}\lesssim K^{1/2-\alpha}$, after which Theorem~\ref{thm:rates} implies $\Phi_b(Q)\lesssim Q^{3/2-\alpha}$.
Together with $\|S_db\|_\infty^2\le\sum_k\|\Delta_kb\|_\infty^2<\infty$, Proposition~\ref{prop:34} below therefore yields $\Phi_b(Q)\lesssim Q^{3/4}$. Hence, for $\frac12<\alpha<\frac34$, the exponent $\frac32-\alpha$ is not sharp.
The sharp exponent for this range remains unknown to us; for chains, however, Theorem~\ref{thmB} yields $Q^{1/2}$. The dyadic version of the sharp extrapolation theorem \cite{DGPP05} transfers bounds of the form $\Phi_b(Q)\lesssim Q^\beta$ to $L^p(w)$, $1<p<\infty$, replacing $Q$ by $[w]_{A_p^d}^{\max(1,1/(p-1))}$. This version is obtained by using the dyadic maximal function in the Rubio de Francia algorithm.
\end{remark}

\section{Nesting and chains}\label{sec:tree}

\begin{proposition}\label{prop:notdelta}
Let $J_k=[2k,2k+2^{-k})$, $I_k=[0,2^{-k})$, $L\ge2$, and
\[
b^{\mathrm{dis}}_L=\sum_{k<L}|J_k|^{1/2}h_{J_k},\qquad b^{\mathrm{ch}}_L=\sum_{k<L}(-1)^k|I_k|^{1/2}h_{I_k}.
\]
For either choice of symbol for $b$, we have $\|\Delta_kb\|_\infty=1$ when $0\le k<L$, and otherwise $\Delta_kb=0$. Nevertheless, for every $Q\ge1$ we have $\Phi_{b^{\mathrm{dis}}_L}(Q)=Q^{1/2}$, whereas $\Phi_{b^{\mathrm{ch}}_L}(L)\ge L/\sqrt8$.
\end{proposition}

\begin{proof}
The pairwise disjoint dyadic intervals $J_k$ yield $\|\pi_{b^{\mathrm{dis}}_L}\|_{L^2(w)}=\max_k(\avg w_{J_k}\avg\sigma_{J_k})^{1/2}\le\Aw^{1/2}$ by Lemma~\ref{lem:single}. Conversely, choose, for a given $Q\ge1$, an $t\ge1$ satisfying $\frac{(t+1)^2}{4t}=Q$; set $w=t$ on the left child of $J_0=[0,1)$, taking $w=1$ everywhere else. Whenever $J$ lies in a child of $J_0$ or does not intersect $J_0$, the function $w$ is constant on $J$. In the remaining case, $J\supseteq J_0$; taking $\theta=|J_0|/(2|J|)\in(0,\frac12]$,
\[
\avg w_J\avg\sigma_J=(1+(t-1)\theta)(1+(t^{-1}-1)\theta)=1+\frac{(t-1)^2}t\,\theta(1-\theta)\le\frac{(t+1)^2}{4t},
\]
with equality attained for $J=J_0$. Consequently, $\Aw=Q$ as well as $\|\pi_{b^{\mathrm{dis}}_L}\|_{L^2(w)}\ge(\avg w_{J_0}\avg\sigma_{J_0})^{1/2}=Q^{1/2}$. The second claim follows by applying Proposition~\ref{prop:chain} with $\eps=1$ and $\delta=1/L$.
\end{proof}

Here we have $\|S_db^{\mathrm{dis}}_L\|_\infty=1$ together with $\|S_db^{\mathrm{ch}}_L\|_\infty=\sqrt L$.

\begin{theorem}\label{thm:chain}
For $b=\sum_kc_k|I_k|^{1/2}h_{I_k}$ along a chain $I_0\supsetneq I_1\supsetneq\cdots$ such that $I_{k+1}$ is a child of $I_k$,
\begin{equation}\label{eq:chain}
\|\pi_b\|_{L^2(u)\to L^2(v)}\le8\,[u,v]_{A_2^d}^{1/2}\Big(\sum_kc_k^2\Big)^{1/2}.
\end{equation}
\end{theorem}

\begin{proof}
By applying the finite-subfamily principle from Section~\ref{sec:prelim}, we need only consider finitely many nonzero coefficients.
 After adding one child to the chain if needed, we may take the chain to be $I_0,\dots,I_N$ with $c_N=0$; thus $I_N$ serves as an auxiliary terminal cell.
 Set $A_j=I_j\setminus I_{j+1}$ for $j<N$ and $A_N=I_N$. The relation $|\avg f_I|\le\avg{|f|}_I$ ensures that the pointwise bound below remains valid upon replacing $a_k$ by $\avg{|f|}_{I_k}$, allowing us to assume $f\ge0$.
 For $k<j$, on $A_j$ the function $h_{I_k}$ takes the value $e_k|I_k|^{-1/2}$, where $e_k\in\{\pm1\}$ does not depend on $j$; on $A_j$, it instead satisfies $|h_{I_j}|=|I_j|^{-1/2}$. Therefore, taking $a_k=\avg f_{I_k}$,

\begin{equation}\label{eq:ptchain}
|\pi_bf|\le|c_j|a_j+\sum_{k<j}|c_k|a_k\quad\text{on }A_j,\qquad\pi_bf=0\ \text{off }I_0 .
\end{equation}
Let $F_i=\int_{A_i}f$, so that $a_k=|I_k|^{-1}\sum_{i\ge k}F_i$; applying the Cauchy--Schwarz inequality then yields $\sum_iF_i^2/\sigma(A_i)\le\|f\|_{L^2(u)}^2$. Separating the sum that defines $a_k$ into the cases $i\ge j$ and $i<j$ and invoking $|I_k|=2^{j-k}|I_j|$, we obtain

\begin{equation}\label{eq:split}
\sum_{k<j}|c_k|a_k=\gamma_ja_j+\sum_{i<j}\frac{\eta_iF_i}{|I_i|},\qquad\gamma_j=\sum_{k<j}|c_k|2^{k-j},\quad\eta_i=\sum_{k\le i}|c_k|2^{k-i},
\end{equation}
and an application of Young's inequality yields $\sum\gamma_j^2\le\sum c_k^2$ together with $\sum\eta_i^2\le4\sum c_k^2$. Vanishing at certain indices is possible for some sequences used below, as permitted by Lemma~\ref{lem:hardy}. To handle $\lambda\in\{|c|,\gamma\}$, use the first inequality of \eqref{eq:hardy}, taking $u_j=\lambda_j^2v(A_j)|I_j|^{-2}$, $v_i=\sigma(A_i)$ and $x_i=F_i$; here $v(A_j)\le v(I_j)\le [u,v]_{A_2^d}|I_j|^2/\sigma(I_j)$ and $\sigma(I_n)\le\sigma(I_j)$ whenever $j\le n$,
\[
B_1=\sup_n\sum_{j\le n}\lambda_j^2\frac{v(A_j)}{|I_j|^2}\sigma(I_n)\le [u,v]_{A_2^d}\sup_n\sum_{j\le n}\lambda_j^2\frac{\sigma(I_n)}{\sigma(I_j)}\le [u,v]_{A_2^d}\sum_j\lambda_j^2,
\]
hence $\sum_jv(A_j)\lambda_j^2a_j^2\le4[u,v]_{A_2^d}\sum\lambda^2\,\|f\|^2$. To treat the final sum, use the second inequality of \eqref{eq:hardy}, taking $u_j=v(A_j)$, $x_i=\eta_iF_i/|I_i|$ and $v_i=\eta_i^2\sigma(A_i)/|I_i|^2$; here $\sum_{j>n}v(A_j)=v(I_{n+1})\le v(I_i)$ whenever $i\le n$ and $\sigma(A_i)\le\sigma(I_i)\le [u,v]_{A_2^d}|I_i|^2/v(I_i)$,
\[
B_2\le [u,v]_{A_2^d}\sup_n\sum_{i\le n}\eta_i^2\frac{v(I_{n+1})}{v(I_i)}\le4[u,v]_{A_2^d}\sum_kc_k^2 .
\]
Substituting \eqref{eq:split} into \eqref{eq:ptchain}, followed by an application of Minkowski's inequality in $\ell^2(v(A_j))$, yields
\[
\|\pi_bf\|_{L^2(v)}\le(2+2+4)\Big([u,v]_{A_2^d}\sum_kc_k^2\Big)^{1/2}\|f\|_{L^2(u)} .
\]
\end{proof}

Along a chain, each $A_j$ carries a constant value of $S_db$, with $S_db^2=\sum_{k\le j}c_k^2$ on that set; consequently, $\|S_db\|_\infty^2=\sum_kc_k^2$. Beyond $I_K$, the chain's tail has $\bmo$ norm bounded above by $\sqrt2\sup_{k>K}|c_k|$; thus, we also have
\begin{equation}\label{eq:chaintail}
\Phi_b(Q)\le\min_K\Big[8Q^{1/2}\Big(\sum_{k\le K}c_k^2\Big)^{1/2}+\sqrt2C_BQ\sup_{k>K}|c_k|\Big].
\end{equation}
In particular, \eqref{eq:chaintail} together with Proposition~\ref{prop:lower} implies $\Phi_b(Q)\simeq Q^{1/2}$ for each nonzero $(c_k)\in\ell^2$. As an example, consider the left-descendant chain symbol $\sum_{k\ge1}\frac1k|I_k|^{1/2}h_{I_k}$ on $[0,1)$: its value on $A_j$ is $\sum_{k<j}\frac1k-\frac1j$, making it unbounded, while it lies in $\cmo$ and exhibits the minimal growth.

We next consider lower bounds along a chain.

\begin{proposition}\label{prop:dualchain}
Let $P=I_0\supsetneq I_1\supsetneq\cdots$ denote an infinite chain such that $I_{k+1}$ is a child of $I_k$, and let $c_k=b_{I_k}|I_k|^{-1/2}$, $Q\ge2$ and $\delta=1/Q$. Then

\begin{equation}\label{eq:dualchain}
\Phi_b(Q)\ge\theta(\delta)\,Q^{1/2}\Big(\sum_{0\le k<Q}c_k^2\Big)^{1/2}-\big(1+\sqrt{3Q}\big)\|b\|_{\bmo},
\end{equation}
where $\theta(\delta)=(\tau_\delta/2)^{1/2}\,r(r-1)/(r+1)$, $r=2^{1-\delta}$ and $\tau_\delta=(1-2^{\delta-2})/(2-\delta)$. Whenever $0<\delta\le\frac12$, we have $\theta(\delta)\ge\frac1{11}$; moreover, $\theta(\delta)\to1/\sqrt{12}$ as $\delta\to0$. The term involving $\|b\|_{\bmo}$ may be omitted if $b$ has no coefficients outside the chain.
\end{proposition}

\begin{proof}
\emph{Reduction to a left chain.} For every $k$ for which $I_{k+1}$ is the right child of $I_k$, interchange the two children of $I_k$ along with all of their dyadic subintervals.
 The intervals $I_k$ contract to a single point in $\overline P$; we disregard this null set. Every remaining point in $P$ undergoes only finitely many exchanges. Consequently, the composition $\tau$ defines a piecewise translation of $P$ that is well-defined, preserves measure, and sends $\D(P)$ onto $\D(P)$.
Set $\tau=\mathrm{id}$ outside $P$. Because $\tau$ transports averages over dyadic subintervals of $P$ while preserving integrals over $P$, it follows that $[w\circ\tau^{-1}]_{A_2^d}=\Aw$. Thus $Uf=f\circ\tau^{-1}$ gives an isometry from $L^2(w)$ onto $L^2(w\circ\tau^{-1})$.
The relations $Uh_I=\pm h_{\tau I}$ for $I\subseteq P$, $Uh_I=h_I$ for every other $I$, and $\avg{Uf}_{\tau I}=\avg f_I$ yield $U\pi_bU^{-1}=\pi_{b'}$ with $|b'_{\tau I}|=|b_I|$. Consequently, $\Phi_{b'}=\Phi_b$ and $\|b'\|_{\bmo}=\|b\|_{\bmo}$. Only $|c_k|$ and $\|b\|_{\bmo}$ enter the statement; we can therefore take the chain to be left-descendant. Following Proposition~\ref{prop:chain}, we present the computation for $P=[0,1)$, $I_k=[0,2^{-k})$.

\emph{Test.} Let $w=x^{\delta-1}\1_P+\1_{\R\setminus P}$. Exchanging the roles of $w$ and $\sigma$ in the argument of Proposition~\ref{prop:chain} gives $\Aw\le1/\delta=Q$. We employ $\|\pi_b\|_{L^2(w)}=\|\pi_b^*\|_{L^2(\sigma)}$ with $\pi_b^*g=\sum_Ib_Ig_I\1_I/|I|$, together with the test function $g=w\1_P$ satisfying $\|g\|_{L^2(\sigma)}^2=w(P)=1/\delta$. Decompose $\pi_b^*g=G_{\mathrm{ch}}+G_{\mathrm{off}}+G_{\mathrm{anc}}$ into the cases where $I=I_k$ for some $k$, where $I\subseteq P$ lies outside the chain, and where $I\supsetneq P$; for intervals disjoint from $P$, we have $g_I=0$.

\emph{The chain.} Because $w(I_m)=2^{-m\delta}/\delta$,
\[
\begin{aligned}
g_{I_k}&=|I_k|^{-1/2}\big(w(I_{k+1})-w(I_k\setminus I_{k+1})\big)\\
&=|I_k|^{-1/2}\big(2w(I_{k+1})-w(I_k)\big)=|I_k|^{-1/2}2^{-k\delta}\frac{r-1}\delta,
\end{aligned}
\]
therefore $G_{\mathrm{ch}}=\frac{r-1}\delta\sum_kc_kr^k\1_{I_k}$. This takes the value $\frac{r-1}\delta r^j(T_rc)_j$ on $A_j=I_j\setminus I_{j+1}$, with $(T_rc)_j=\sum_{k\le j}c_kr^{k-j}$. The relation $c_j=(T_rc)_j-r^{-1}(T_rc)_{j-1}$ implies $\|(c_k)_{k<J}\|_{\ell^2}\le(1+r^{-1})\|((T_rc)_j)_{j<J}\|_{\ell^2}$ for each $J$. In addition, $\sigma(A_j)=\tau_\delta2^{-j(2-\delta)}$, while $r^{2j}2^{-j(2-\delta)}=2^{-j\delta}>\frac12$ for $j<J:=\lceil Q\rceil$. Thus
\begin{equation}\label{eq:Gch}
\|G_{\mathrm{ch}}\|_{L^2(\sigma)}^2\ge\frac{(r-1)^2}{\delta^2}\,\frac{\tau_\delta}2\sum_{j<J}(T_rc)_j^2\ge\theta(\delta)^2\,Q\Big(\sum_{k<Q}c_k^2\Big)\|g\|_{L^2(\sigma)}^2 .
\end{equation}

\emph{Off the chain.} Every $I\subseteq P$ not belonging to the chain belongs to some $A_j$. The relation $|w'|\le2w_j/|A_j|$ holds on $A_j$ with $w_j=\min_{A_j}w$, yielding $|g_I|\le\frac14|I|^{3/2}\sup_I|w'|\le|I|^{3/2}w_j/(2|A_j|)$. Combining $|b_I|\le\|b\|_{\bmo}|I|^{1/2}$ with $\sum_{I\subseteq A_j,\,I\ni x}|I|\le2|A_j|$ yields $|G_{\mathrm{off}}|\le\|b\|_{\bmo}w_j$ on $A_j$; moreover, $w_j^2\sigma\le w$ on $A_j$ implies $\|G_{\mathrm{off}}\|_{L^2(\sigma)}\le\|b\|_{\bmo}\|g\|_{L^2(\sigma)}$.

\emph{Ancestors.} For $A\supsetneq P$, the support of $g$ lies in one child of $A$, which gives $|g_A|=|A|^{-1/2}/\delta$; consequently
\[
|G_{\mathrm{anc}}(x)|\le\|b\|_{\bmo}\,\delta^{-1}\sum_{A\supsetneq P,\,A\ni x}|A|^{-1}.
\]
This has upper bound $\|b\|_{\bmo}/\delta$ on $P$, with $\sigma(P)\le1$. Taking $A^{(m)}$ to be the ancestor of length $2^m$ gives the upper bound $\|b\|_{\bmo}2^{1-m}/\delta$ on $A^{(m)}\setminus A^{(m-1)}$; this set has measure $2^{m-1}$ and satisfies $\sigma=1$. Summation yields,
\begin{equation}\label{eq:Ganc}
\|G_{\mathrm{anc}}\|_{L^2(\sigma)}^2\le\frac{3\|b\|_{\bmo}^2}{\delta^2}=3Q\|b\|_{\bmo}^2\|g\|_{L^2(\sigma)}^2 .
\end{equation}

The bound for $G_{\mathrm{off}}$, together with \eqref{eq:Ganc}, \eqref{eq:Gch} and the triangle inequality, yields \eqref{eq:dualchain}. When all coefficients of $b$ off the chain vanish, $G_{\mathrm{off}}=G_{\mathrm{anc}}=0$ follows. Finally, given $0<\delta\le\frac12$, we have $r\in[\sqrt2,2)$ and $\tau_\delta\ge\frac12(1-2^{-3/2})$, while the function $r\mapsto r(r-1)/(r+1)$ is increasing; therefore $\theta(\delta)\ge(\frac14(1-2^{-3/2}))^{1/2}\frac{\sqrt2(\sqrt2-1)}{\sqrt2+1}>\frac1{11}$. Because $\delta\to0$, $r\to2$ and $\tau_\delta\to\frac38$, we conclude $\theta(\delta)\to(\frac3{16})^{1/2}\cdot\frac23=\frac1{\sqrt{12}}$.
\end{proof}

\begin{proof}[Proof of Theorem~\ref{thmB}]
(i) Theorem~\ref{thm:chain} supplies the upper bound; $(\sum c_k^2)^{1/2}=\|S_db\|_\infty$ has already been observed. To obtain the lower bound, first extend any finite chain using zero coefficients, then invoke \eqref{eq:dualchain} with the error term omitted: $[b]_{\G_{1/2}}\ge Q^{-1/2}\Phi_b(Q)\ge\theta(1/Q)(\sum_{k<Q}c_k^2)^{1/2}$, and let $Q\to\infty$.

(ii) We can suppose that $[b]_{\G_{1/2}}<\infty$, giving $\|b\|_{\bmo}\le\Phi_b(1)\le[b]_{\G_{1/2}}$. Choose $x\in\R$ and $P\in\D$ containing $x$. An infinite chain is formed by the dyadic intervals $I\subseteq P$ containing $x$, so Proposition~\ref{prop:dualchain} yields
\[
\begin{aligned}
\theta(1/Q)\Big(\sum_{k<Q}c_k^2\Big)^{1/2}&\le Q^{-1/2}\Phi_b(Q)+(Q^{-1/2}+\sqrt3)\|b\|_{\bmo}\\
&\le\big(1+Q^{-1/2}+\sqrt3\big)[b]_{\G_{1/2}} .
\end{aligned}
\]
For fixed $P$, put $Q\to\infty$, thereby bounding $(\sum_{I\subseteq P,\,I\ni x}b_I^2/|I|)^{1/2}$ by $(1+\sqrt3)\sqrt{12}\,[b]_{\G_{1/2}}$.
Next, allow $P$ to increase through its ancestors. Monotone convergence then yields $S_db(x)\le(1+\sqrt3)\sqrt{12}\,[b]_{\G_{1/2}}$; thus (ii) holds with $C=(1+\sqrt3)\sqrt{12}<10$.
\end{proof}

For chains whose coefficients are monotone, the two estimates have matching growth rates, determining the growth exactly.

\begin{theorem}\label{thm:monotone}
Let $I_0\supsetneq I_1\supsetneq\cdots$, where $I_{k+1}$ is a child of $I_k$; take $a_0\ge a_1\ge\cdots\ge0$ with $a_k\to0$ and $\eps_k\in\{-1,1\}$, and
\[
b=\sum_{k\ge0}\eps_ka_k|I_k|^{1/2}h_{I_k},\qquad A_N=\Big(\sum_{k<N}a_k^2\Big)^{1/2}.
\]
Then the tails $b^{(N)}=\sum_{k\ge N}\eps_ka_k|I_k|^{1/2}h_{I_k}$ obey $\|b^{(N)}\|_{\bmo}\le\sqrt2\,a_N$, and $b\in\cmo$; consequently, $\|b^{>K}\|_{\bmo}\le\sqrt2\,a_{K+1}$ for $K\ge0$ provided $|I_0|=1$, while for $Q\ge2$
\begin{equation}\label{eq:monotone}
c\,Q^{1/2}A_{\lceil Q\rceil}\le\Phi_b(Q)\le C\,Q^{1/2}A_{\lceil Q\rceil}
\end{equation}
with absolute constants $c,C>0$. Under the assumptions $I_k=[0,2^{-k})$ and $\eps_k=(-1)^k$, we also have $\|b\|_\infty\le2a_0$.
\end{theorem}

\begin{proof}
Only chain intervals have coefficients. Thus, for $N\ge0$, the Carleson sum of $b^{(N)}$ over $J\in\D$ is zero except when $J=I_m$ for some $m\ge0$ or $J\supsetneq I_0$; in the latter situation, it does not exceed the sum for $J=I_0$. When $J=I_m$ and $s=\max(N,m)$, monotonicity yields
\[
\frac1{|I_m|}\sum_{k\ge s}a_k^2|I_k|=\sum_{k\ge s}a_k^2\,2^{m-k}\le a_s^2\,2^{m-s+1}\le2a_N^2,
\]
hence $\|b^{(N)}\|_{\bmo}\le\sqrt2\,a_N$. Under the condition $|I_0|=1$, we obtain $b^{>K}=b^{(K+1)}$. As $b-b^{(N)}$ is a finite Haar sum with $a_N\to0$, it follows that $b\in\cmo$.

For the upper bound, take $N=\lceil Q\rceil$ and $\Aw\le Q$. Applying estimate \eqref{eq:chain} to $b-b^{(N)}$, together with \eqref{eq:bez} to $b^{(N)}$, yields
\begin{equation}\label{eq:monoupper}
\|\pi_b\|_{L^2(w)}\le8Q^{1/2}A_N+\sqrt2\,C_BQ\,a_N .
\end{equation}
By monotonicity we have $A_N^2\ge Na_{N-1}^2\ge Qa_N^2$ and therefore $Qa_N\le Q^{1/2}A_N$. Hence \eqref{eq:monoupper} yields the upper bound in \eqref{eq:monotone} with $C=8+\sqrt2\,C_B$. To obtain the lower bound, use \eqref{eq:dualchain} with the error term omitted: $\Phi_b(Q)\ge\theta(1/Q)Q^{1/2}(\sum_{0\le k<Q}a_k^2)^{1/2}$, where $\theta(1/Q)\ge\frac1{11}$ and $\{k\ge0:k<Q\}=\{0,\dots,\lceil Q\rceil-1\}$. Thus the choice $c=\frac1{11}$ suffices.

Finally, assuming $I_k=[0,2^{-k})$, on $A_j=I_j\setminus I_{j+1}$ we obtain $|I_k|^{1/2}h_{I_k}=1$ for $k<j$ as well as $=-1$ for $k=j$; consequently, $b=\sum_{k<j}(-1)^ka_k-(-1)^ja_j$ there. For an alternating series whose terms are nonincreasing, the partial sums belong to $[0,a_0]$, which gives $|b|\le2a_0$.
\end{proof}

\begin{corollary}\label{cor:critical}
For $\kappa>0$ take $b_\kappa=\sum_{k\ge0}(-1)^k(1+k)^{-\kappa}|I_k|^{1/2}h_{I_k}$ subject to $I_k=[0,2^{-k})$. We then have $b_\kappa\in\cmo\cap L^\infty$ and, for $K\ge0$, $\|b_\kappa^{>K}\|_{\bmo}\le\sqrt2(1+K)^{-\kappa}$, whereas for $Q\ge2$
\begin{equation}\label{eq:kappa}
\Phi_{b_\kappa}(Q)\simeq_\kappa\begin{cases}Q^{1-\kappa},&0<\kappa<\frac12,\\(Q\log(e+Q))^{1/2},&\kappa=\frac12,\\Q^{1/2},&\kappa>\frac12.\end{cases}
\end{equation}
In particular, replacing $Q^{1/2}\log(e+Q)$ by $Q^{1/2}$ in the case $\kappa=\frac12$ of Theorem~\ref{thmA}\textup{(ii)} is impossible, even for bounded symbols in $\cmo$; precisely when $\alpha>\frac12$, we have $b_{1/2}\in\G_\alpha$, and $S_db_{1/2}\notin L^\infty$.
\end{corollary}

\begin{proof}
Use Theorem~\ref{thm:monotone} with $a_k=(1+k)^{-\kappa}$ to obtain $\sqrt2\,a_{K+1}\le\sqrt2(1+K)^{-\kappa}$ and $\|b_\kappa\|_\infty\le2$. Integral comparison yields, for $\kappa<\frac12$, $A_N^2\simeq_\kappa N^{1-2\kappa}$; for $\kappa=\frac12$, $\log(1+N)\le A_N^2\le1+\log N$; and for $\kappa>\frac12$, $1\le A_N^2\le\zeta(2\kappa)$. Now substitute $N=\lceil Q\rceil$. When $\kappa=\frac12$, one has $\|S_db_{1/2}\|_\infty^2=\sum_k(1+k)^{-1}=\infty$, in agreement with Theorem~\ref{thmB}(ii).
\end{proof}

For $0<\kappa<\frac12$, the sharpness asserted in Theorem~\ref{thmA}(ii) is realized by a single bounded chain, as \eqref{eq:kappa} demonstrates. The same exponent, together with explicit constants, follows from Proposition~\ref{prop:rate-sharp}.

\begin{remark}[Partial upper bound]\label{rem:partial}
Denote by $C_{\mathrm{LSF}}(w)$ the best constant for the lower square function estimate $\|F\|_{L^2(w)}\le C_{\mathrm{LSF}}(w)(\sum_IF_I^2\avg w_I)^{1/2}$ on finite Haar series $F$. The square function expression associated with $F=\pi_bf$ is $\sum_I(b_I\avg f_I)^2\avg w_I$. Thus, upon writing $f=\varphi\sigma$ and using the weighted dyadic Carleson embedding theorem for $\alpha_I=b_I^2\avg w_I\avg\sigma_I^2$, with Carleson constant relative to $\sigma$ bounded above by
\begin{equation}\label{eq:carl}
\sup_J\frac{\Aw}{\sigma(J)}\sum_{I\subseteq J}b_I^2\avg\sigma_I=\sup_J\frac{\Aw}{\sigma(J)}\int_J\sigma\sum_{I\subseteq J}\frac{b_I^2}{|I|}\1_I\le\Aw\|S_db\|_\infty^2,
\end{equation}
it follows from \eqref{eq:carl} that $\|\pi_b\|_{L^2(w)}\le2C_{\mathrm{LSF}}(w)\Aw^{1/2}\|S_db\|_\infty$.
In the case of the standard dyadic filtration, $C_{\mathrm{LSF}}(w)\le C[w]_{A_\infty}^{1/2}$; here $[w]_{A_\infty}$ denotes the dyadic Fujii--Wilson constant \eqref{eq:Ainf}.
The martingale characteristic $[w]_{\infty,\mathcal D}$ of \cite[Definition~2.6]{DIPTV} agrees with this constant. For homogeneous filtrations, \cite[Remark~3.5]{DIPTV} states the estimate in this form, whereas \cite[Theorem~3.4]{DIPTV} supplies it using the larger semiclassical characteristic.
In the standard dyadic filtration, this is also the dyadic version of \cite[Lemma~2.3]{Wi89} for $p=2$ and $\eta=1$.
Consequently, \eqref{eq:Sq} holds up to a factor $[w]_{A_\infty}^{1/2}$. If only the weaker hypothesis $b\in\bmo$ is imposed, the Carleson step incurs an extra factor $[\sigma]_{A_\infty}$.
For the weights $w_\delta$ in Proposition~\ref{prop:chain}, the bound is sharp, and their $A_\infty$ characteristic remains bounded uniformly in $\delta$. For the weights in Proposition~\ref{prop:dualchain}, however, a factor $Q^{1/2}$ is lost: Theorem~\ref{thm:chain} identifies $Q^{1/2}$ as the correct order on chains for these weights. This loss is removed by Corollary~\ref{cor:power} below. In particular, when $0\le\gamma\le1$,
\[
\|\pi_b\|_{L^2(w)}\le C\,\Aw^{(1+\gamma)/2}\,\|S_db\|_\infty\qquad\text{whenever }[w]_{A_\infty}\le\Aw^{\gamma},
\]
which, for improves the exponents in \eqref{eq:bez}, at the cost of replace $\|b\|_{BMO^d}$ by the larger quantity $\|S_d b\|_{\infty}$
\end{remark}

We end this section with an examination of the hypothesis of \eqref{eq:Sq}, a proof of \eqref{eq:Sq} for every symbol on a class of weights containing the power weights used in our lower bounds, and a proof of a sublinear bound for all weights.
This hypothesis takes the form of a Carleson condition uniform over all measures.

\begin{lemma}\label{lem:univ}
For each positive locally finite Borel measure $\mu$ on $\R$ and each $J\in\D$,
\begin{equation}\label{eq:univ}
\sum_{I\subseteq J}\frac{b_I^2}{|I|}\,\mu(I)\le\|S_db\|_\infty^2\,\mu(J),
\end{equation}
and the smallest constant for which this holds is $\|S_db\|_\infty^2$.
\end{lemma}

\begin{proof}
The expression on the left is $\int_J\sum_{I\subseteq J}(b_I^2/|I|)\1_I\,d\mu\le\int_J(S_db)^2\,d\mu$. For the converse, suppose $\lambda<\|S_db\|_\infty^2$.
With $c_I=b_I|I|^{-1/2}$, the set of points $x$ where $\sum_{I\ni x}c_I^2>\lambda$ then has positive measure.
Choose such an $x$ away from the dyadic endpoints, and take $\mu$ to be the point mass at $x$. For $J$ a sufficiently large ancestor of $x$, the left side of \eqref{eq:univ} is greater than $\lambda\mu(J)$.
\end{proof}

Upon taking $d\mu=\sigma\,dx$ and including the factor $\Aw$, \eqref{eq:univ} becomes the Carleson condition \eqref{eq:carl}. Observe also that $\{b:S_db\in L^\infty\}$ is a solid class: for each $\mathcal F\subseteq\D$, $S_d(b_{\mathcal F})\le S_db$ holds pointwise, whereas restriction of coefficients on $L^\infty$ is unbounded. Unbounded functions belong to this class, including the chain symbol with $c_k=1/k$ following Theorem~\ref{thm:chain}; some bounded functions do not, including $b_{1/2}$ from Corollary~\ref{cor:critical}.

For the weights, we keep track of the characteristic's local imbalance. Set
\begin{equation}\label{eq:imb}
d_I(w)=\frac{w(I_l)-w(I_r)}{w(I)},\qquad e_I=d_I(w)+d_I(\sigma),
\end{equation}
so that $|d_I(w)|<1$ as well as $\avg w_{I_l}=\avg w_I(1+d_I(w))$, $\avg w_{I_r}=\avg w_I(1-d_I(w))$, with the same holding for $\sigma$. Taking the product of the averages, with $q_I=\avg w_I\avg\sigma_I$,
\begin{equation}\label{eq:qchild}
q_{I_l}=q_I\big(1+d_I(w)\big)\big(1+d_I(\sigma)\big),\qquad q_{I_r}=q_I\big(1-d_I(w)\big)\big(1-d_I(\sigma)\big),\qquad e_I=\frac{q_{I_l}-q_{I_r}}{2q_I},
\end{equation}
so the vanishing of $e_I$ occurs precisely when the characteristic is the same on both children of $I$.

\begin{proposition}\label{prop:drift}
Let $\widehat h_I=h_I-d_I(w)\,|I|^{-1/2}\1_I$. In $L^2(w)$, the functions $\widehat h_I$ are pairwise orthogonal with $\|\widehat h_I\|_{L^2(w)}^2=(1-d_I(w)^2)\avg w_I$, and
\begin{equation}\label{eq:dec}
\pi_bf=U_bf+D_bf,\qquad U_bf=\sum_Ib_I\avg f_I\,\widehat h_I,\qquad D_bf=\sum_Ic_I\,d_I(w)\avg f_I\1_I,
\end{equation}
where both operators depend on $w$ and $c_I=b_I|I|^{-1/2}$. In addition, $\|U_b\|_{L^2(w)}\le2\,\Aw^{1/2}\|S_db\|_\infty$.
\end{proposition}

\begin{proof}
From $\int_Ih_Iw=|I|^{-1/2}(w(I_l)-w(I_r))=d_I(w)|I|^{-1/2}w(I)$ we obtain $\int\widehat h_I\,w=0$, while $b_I\,d_I(w)|I|^{-1/2}=c_I\,d_I(w)$ yields \eqref{eq:dec}. When $J\subsetneq I$ holds, constancy of $\widehat h_I$ on $J$ implies $\int\widehat h_I\widehat h_J\,w=0$; the supports associated with disjoint intervals are disjoint. Expansion gives $\|\widehat h_I\|_{L^2(w)}^2=\avg w_I-2d_I(w)^2\avg w_I+d_I(w)^2\avg w_I$. Thus, taking $f=\varphi\sigma$ and $\avg\varphi_{I,\sigma}=\sigma(I)^{-1}\int_I\varphi\sigma$,
\[
\|U_bf\|_{L^2(w)}^2\le\sum_Ib_I^2\avg f_I^2\avg w_I=\sum_Ib_I^2\avg w_I\avg\sigma_I^2\,\avg\varphi_{I,\sigma}^2\le4\,\Aw\,\|S_db\|_\infty^2\,\|f\|_{L^2(w)}^2
\]
using \eqref{eq:carl} together with the weighted Carleson embedding theorem, as in Remark~\ref{rem:partial}.
\end{proof}

\begin{theorem}\label{thm:spine}
Let $w\in A_2^d$, and take a chain $\mathcal C\subseteq\D$, namely a totally ordered family that includes every dyadic interval containing any member of the family. If
\begin{equation}\label{eq:energy}
E=\operatorname*{ess\,sup}_{x\in\R}\ \sum_{I\ni x,\ I\notin\mathcal C}e_I^2<\infty,
\end{equation}
then, for each symbol $b$,
\[
\|\pi_b\|_{L^2(w)}\le C\,e^{E}\,\Aw^{1/2}\,\|S_db\|_\infty ,
\]
where $C$ is an absolute constant.
\end{theorem}

\begin{proof}
Write $Q=\Aw$ and $M=\|S_db\|_\infty$, and decompose the coefficients as $b=b_{\mathcal C}+b'$ with $b_{\mathcal C}=\sum_{I\in\mathcal C}b_Ih_I$. Solidity gives $S_d\le M$ for both parts. To see that the chain is consecutive, take $J\subsetneq I$ in $\mathcal C$: the child of $I$ that contains $J$ is an ancestor of $J$ and consequently belongs to $\mathcal C$. Thus any finite subfamily of $\mathcal C$ is contained in a single descending chain of children. Its smallest member is contained in every interval in that subfamily, giving $\sum_{I\in\mathcal C}c_I^2\le M^2$. Applying Theorem~\ref{thm:chain} together with the finite-subfamily principle from Section~\ref{sec:prelim} yields $\|\pi_{b_{\mathcal C}}\|_{L^2(w)}\le8\,Q^{1/2}M$.

Proposition~\ref{prop:drift} gives $\pi_{b'}=U_{b'}+D_{b'}$ as well as $\|U_{b'}\|_{L^2(w)}\le2Q^{1/2}M$. Taking $f=\varphi\sigma$ yields the pointwise bound $|D_{b'}f|\le T(|\varphi|\sigma)$, where
\[
T(\varphi\sigma)=\sum_{I\notin\mathcal C}\lambda_I\avg{\varphi\sigma}_I\1_I,\qquad\lambda_I=|c_I\,d_I(w)| .
\]
We check the testing conditions in Theorem~\ref{thm:LSU} for $T$. Take $V_I=w(I)\avg\sigma_I^2$ and $\widetilde V_I=\sigma(I)\avg w_I^2$. Expansion of $V_{I_\pm}=\tfrac12w(I)(1\pm d_I(w))\avg\sigma_I^2(1\pm d_I(\sigma))^2$ and of its dual gives
\begin{equation}\label{eq:potid}
V_{I_l}+V_{I_r}=\big(1-d_I(w)^2+e_I^2\big)V_I,\qquad
\widetilde V_{I_l}+\widetilde V_{I_r}=\big(1-d_I(w)^2+2d_I(w)e_I\big)\widetilde V_I .
\end{equation}
Let $R\notin\mathcal C$. Because $\mathcal C$ includes all ancestors of its members, no $I\subseteq R$ belongs to $\mathcal C$. Given $I\subseteq R$, set $D_I=\prod_{R\supseteq A\supsetneq I}(1+e_A^2)$; consequently $D_R=1$, and both children $K$ of $I$ satisfy $D_K=D_I(1+e_I^2)$. Moreover, \eqref{eq:energy} gives $D_I(1+e_I^2)\le e^{E}$, since every interval in the product contains $I$ and lies off $\mathcal C$. Using the first identity of \eqref{eq:potid},
\[
\frac{V_I}{D_I}-\frac{V_{I_l}+V_{I_r}}{D_I(1+e_I^2)}=\frac{V_I}{D_I}\,\frac{d_I(w)^2}{1+e_I^2}\ge e^{-E}\,d_I(w)^2\,V_I ,
\]
and the sum over $I\subseteq R$ telescopes, yielding an upper bound of $V_R$. The computation also applies to $\widetilde V$ upon using $\widetilde D_I=\prod(1+2e_A^2)\le e^{2E}$ and $1-d^2+2de\le1-\tfrac12d^2+2e^2$. Thus
\begin{equation}\label{eq:teles}
\sum_{I\subseteq R}d_I(w)^2\,V_I\le e^{E}\,V_R,\qquad
\sum_{I\subseteq R}d_I(w)^2\,\widetilde V_I\le2e^{2E}\,\widetilde V_R .
\end{equation}
For $J\in\D$, pairwise disjointness holds for the maximal intervals of $\{I\subseteq J:I\notin\mathcal C\}$; moreover, $V_R\le Q\sigma(R)$ and $\widetilde V_R\le Qw(R)$. Consequently, \eqref{eq:teles} yields
\[
\sum_{I\subseteq J,\,I\notin\mathcal C}d_I(w)^2V_I\le e^{E}Q\,\sigma(J),\qquad
\sum_{I\subseteq J,\,I\notin\mathcal C}d_I(w)^2\widetilde V_I\le2e^{2E}Q\,w(J).
\]
Using $\sum_{I\ni x}c_I^2\le M^2$ and applying the Cauchy--Schwarz inequality on every path,
\[
\int_J\Big(\sum_{I\subseteq J,\,I\notin\mathcal C}\lambda_I\avg\sigma_I\1_I\Big)^2w
\le M^2\sum_{I\subseteq J,\,I\notin\mathcal C}d_I(w)^2\avg\sigma_I^2\,w(I)\le e^{E}QM^2\sigma(J),
\]
so $\mathfrak T\le e^{E/2}Q^{1/2}M$; replacing $\avg\sigma_I$ and $w$ by $\avg w_I$ and $\sigma$ in the same argument yields $\mathfrak T^*\le\sqrt2\,e^{E}Q^{1/2}M$. By Theorem~\ref{thm:LSU}, we have $\|D_{b'}\|_{L^2(w)}\le C_T(e^{E/2}+\sqrt2e^E)Q^{1/2}M$. Together, the three pieces are bounded above by $(10+3C_T)e^{E}Q^{1/2}M$.
\end{proof}

With $\mathcal C=\emptyset$, weights whose children are balanced or nearly balanced everywhere fall within Theorem~\ref{thm:spine}, regardless of how strongly they branch. One direction of persistent imbalance is absorbed by the chain; this is precisely the configuration appearing in our lower bounds.

\begin{corollary}\label{cor:power}
Let $P\in\D$ and $0<\delta<1$, and set $w=((x-\inf P)/|P|)^{\pm(1-\delta)}$ on $P$ and $w=1$ off $P$, or take its reflection. Then
\[
\|\pi_b\|_{L^2(w)}\le C\,\Aw^{1/2}\,\|S_db\|_\infty
\]
for every $b$, where $C$ is an absolute constant. In particular, at the weights in Propositions~\ref{prop:chain} and~\ref{prop:dualchain}, every symbol satisfies \eqref{eq:Sq}.
\end{corollary}

\begin{proof}
By affine invariance we may choose $P=[0,1)$, placing the singular point at $0$; thus $w=x^{\gamma}$ on $P$, where $\gamma=\pm(1-\delta)$. Define $\mathcal C$ to contain $I_k=[0,2^{-k})$, $k\ge0$, and all dyadic ancestors of $P$. In the terminology of Theorem~\ref{thm:spine}, this forms a chain. If $I\notin\mathcal C$ and $e_I\ne0$, we obtain $I\subseteq A_j=I_j\setminus I_{j+1}$ for some $j\ge0$, since intervals disjoint from $P$ carry $w\equiv1$, whereas any interval meeting both $P$ and its complement is an ancestor of $P$. For $0<h\le a$, put $\rho=h/a$. On $[a,a+h)$, the characteristic of $x^{\gamma}$ is
\[
q_{[a,a+h)}=\Phi_\delta(\rho)=\frac{(1+\rho)^{\delta}-1}{\delta\rho}\cdot\frac{(1+\rho)^{2-\delta}-1}{(2-\delta)\rho}
\]
for either sign of $\gamma$, because interchanging $w$ and $\sigma$ leaves $q_I$ unchanged. As a function of $(\rho,\delta)\in[0,1]\times[0,1]$, $\Phi$ admits a $C^1$ extension with $\Phi\ge1$; consequently, $C_0=\sup|\partial_\rho\Phi_\delta(\rho)|$ is finite and absolute. The parameters of the children of $[a,a+h)$ are $\rho/2$ and $\rho/(2+\rho)$. The last identity in \eqref{eq:qchild} therefore gives
\[
|e_I|\le\frac{C_0}2\Big(\frac\rho2-\frac\rho{2+\rho}\Big)=\frac{C_0\,\rho^2}{4(2+\rho)}\le\frac{C_0\rho^2}8 .
\]
If $x\in A_j$, all intervals $I\ni x$ off $\mathcal C$ are contained in $A_j$. The interval of length $2^{-m}|A_j|$ satisfies $a\ge|A_j|$, which implies $\rho\le2^{-m}$. Thus
\[
\sum_{I\ni x,\,I\notin\mathcal C}e_I^2\le\frac{C_0^2}{64}\sum_{m\ge0}16^{-m}\le\frac{C_0^2}{60},
\]
and with $E\le C_0^2/60$, Theorem~\ref{thm:spine} is applicable.
\end{proof}

The chain can instead be replaced by a family of chains whose tops are disjoint, thereby covering weights that have many separated singular points.

\begin{theorem}\label{thm:forest}
Let $w\in A_2^d$, and suppose $\mathcal C=\bigcup_i\mathcal C_i$ is a union of consecutive chains with pairwise disjoint maximal intervals $T_i$. Here each chain is a nested sequence of dyadic intervals, with each interval followed by one of its children. For $\mathcal C$, denote the quantity \eqref{eq:energy} by $E$. Then, for every $b$,
\[
\|\pi_b\|_{L^2(w)}\le C\,e^{E}\,\Aw^{1/2}\,\|S_db\|_\infty ,
\]
where the constant $C$ is absolute.
\end{theorem}

\begin{proof}
The proof of Theorem~\ref{thm:spine} carries over after two modifications. First, since $\pi_{b_{\mathcal C_i}}f$ depends on $f\1_{T_i}$ and has support in $T_i$, the operators $\pi_{b_{\mathcal C_i}}$ map between orthogonal subspaces of $L^2(w)$ and $\|\pi_{b_{\mathcal C}}\|\le\sup_i\|\pi_{b_{\mathcal C_i}}\|\le8\Aw^{1/2}M$; every chain lies on a single path.

Second, chains in the testing estimates are not required to contain their ancestors. Accordingly, decompose the family $\{I\subseteq J:I\notin\mathcal C\}$ into components for the parent--child relation. In any such component, each interval and all its ancestors through the component root lie off $\mathcal C$; thus the telescoping \eqref{eq:teles} can be applied component by component.
The roots consist of $J$ itself along with those intervals having a parent in $\mathcal C$. There cannot be two nested roots $R'\subsetneq R$ with parents $p'\in\mathcal C_j$, $p\in\mathcal C_k$. Indeed, $p'\subsetneq R\subseteq p$ implies $T_j\cap T_k\ne\emptyset$ and hence $j=k$. The consecutive chain running down from $p$ to $p'$ therefore meets the ancestor of $p'$ with size $|R|$, namely $R$, in contradiction to $R\notin\mathcal C$.
Consequently, the roots whose parent belongs to $\mathcal C$ constitute an antichain. Together with $J$, they decompose into two disjoint layers, giving $\sum_RV_R\le2\Aw\,\sigma(J)$ and $\sum_R\widetilde V_R\le2\Aw\,w(J)$.
The testing constants are therefore doubled, $\mathfrak T\le\sqrt2e^{E/2}\Aw^{1/2}M$ and $\mathfrak T^*\le2e^{E}\Aw^{1/2}M$, yielding the bound with $C=10+4C_T$.
\end{proof}

\begin{remark}\label{rem:blocks}
Theorem~\ref{thm:forest} is applicable to the weight $w$ from the proof of Theorem~\ref{thm:ess}\textup{(ii)}, with $E$ bounded by an absolute constant. This weight agrees with $1$ outside the blocks $P_n=[10n,10n+1)$, $n\in\Z$, and on each such block has the profile of Proposition~\ref{prop:chain} with $\delta=1/Q$. Take $\mathcal C_n$ to be the left-descendant chain of $P_n$ and $\mathcal C=\bigcup_n\mathcal C_n$; its tops $P_n$ are disjoint.
For a fixed point $x$, we estimate the sum in \eqref{eq:energy} over the dyadic intervals that contain it.

Under the conditions $|I|<1$ and $I\notin\mathcal C$, either $w$ is constant on $I$, giving $e_I=0$, or $I$ is contained in an annulus of the block that contains $x$. As in the proof of Corollary~\ref{cor:power}, the contribution of these intervals is at most $C_0^2/60$. When $|I|=1$, either $I$ is a block and thus belongs to $\mathcal C$, or $w\equiv1$ on $I$. Choose $|I|=2^\ell$ with $\ell\ge1$; then $I$ and both of its halves are unions of unit intervals. Denote the numbers of blocks in the two halves of $I$ by $N_l$ and $N_r$, and let $N=N_l+N_r$. Both halves are intervals of $2^{\ell-1}$ consecutive integers. Among $L$ consecutive integers, the number of multiples of $10$ is either $\lfloor L/10\rfloor$ or $\lceil L/10\rceil$. Consequently, $|N_l-N_r|\le1$, with $N\ge\lfloor2^\ell/10\rfloor\ge2^\ell/20$ whenever $\ell\ge4$. On each block the $w$-mass is $1-a_\delta$, where $a_\delta=(1-\delta)/(2-\delta)\in[0,\frac12]$, and the $\sigma$-mass is $1/\delta=Q$; on every remaining unit interval, both masses are $1$.
It follows that
\[
|d_I(w)|=\frac{|N_l-N_r|\,a_\delta}{2^\ell-Na_\delta}\le2^{-\ell},\qquad
|d_I(\sigma)|=\frac{|N_l-N_r|\,(Q-1)}{2^\ell+N(Q-1)}\le\frac{|N_l-N_r|}{\max(N,1)},
\]
thus $|e_I|\le21\cdot2^{-\ell}$ when $\ell\ge4$, whereas $|e_I|<2$ when $\ell=1,2,3$. The sum over $\ell$ therefore yields $E\le12+441\sum_{\ell\ge4}4^{-\ell}+C_0^2/60$, independently of $Q$. Applying Theorem~\ref{thm:forest} now gives $\|\pi_b\|_{L^2(w)}\le C\Aw^{1/2}\|S_db\|_\infty$ at these weights for every $b$.

The step weights in Propositions~\ref{prop:lower} and~\ref{prop:notdelta} satisfy $e_I\ne0$ only on a single dyadic interval together with its ancestors; these constitute a chain as defined in Theorem~\ref{thm:spine}. For the weights $w_t$ in Theorem~\ref{thm:ess}\textup{(iii)}, $e_I\ne0$ occurs only on the unit intervals, which constitute one-element chains with disjoint tops. In each case, $E=0$.
Combining this with Corollary~\ref{cor:power}, which applies to the weights in Propositions~\ref{prop:chain}, \ref{prop:rate-sharp} and~\ref{prop:dualchain} as well as Theorem~\ref{thm:norate}, establishes \eqref{eq:Sq} for every symbol at every weight appearing in this paper's lower bounds.
\end{remark}

Finally, the decomposition in Proposition~\ref{prop:drift} yields a sublinear bound for general weights. Set $\nu_I=(w(I)\sigma(I))^{1/2}=q_I^{1/2}|I|$, write $d_I=d_I(w)$, and take $q_I=\avg w_I\avg\sigma_I$ as in \eqref{eq:qchild}.
Four energy estimates form the basis of the proof; each follows by telescoping an appropriate potential over the dyadic tree.
For the least obvious estimate, Lemma~\ref{lem:energy}(ii), the potential is $(21\Aw-q_I)\nu_I$. The following inequality, applied to the two children of an interval, controls its one-step decrease.

\begin{lemma}\label{lem:twopoint}
Let $x=(x_l,x_r)$ and $u=(u_l,u_r)$ be unit vectors with nonnegative entries in $\R^2$, and suppose $h=1-x\cdot u$ and $d=x_l^2-x_r^2$. We then have
\begin{equation}\label{eq:bellloc}
21h+4\sum_{j\in\{l,r\}}(x_ju_j)^3-1\ge d^2 .
\end{equation}
\end{lemma}

\begin{proof}
Convexity of $t\mapsto t^3$ gives $(x_ju_j)^3\ge x_j^6+3x_j^5(u_j-x_j)$; taking $p=x_l^2$ yields $4(p^3+(1-p)^3)-1=3(2p-1)^2=3d^2$. Consequently, with $z=(x_l^5,x_r^5)$, $4\sum_j(x_ju_j)^3-1\ge3d^2+12\,z\cdot(u-x)$. From $|x|=|u|=1$ it follows that $h\ge0$, $(u-x)\cdot x=-h$ and $|u-x|^2=2h$. Along $x$, the component of $z$ is $(\sum_jx_j^6)\,x$, where $0\le\sum_jx_j^6\le1$; orthogonal to $x$, its component has length $|z\cdot(x_r,-x_l)|=x_lx_r|x_l^4-x_r^4|=x_lx_r|d|\le|d|/2$. It follows that $z\cdot(u-x)\ge-h-\frac12|d|(2h)^{1/2}$, and
\[
21h+4\sum_j(x_ju_j)^3-1\ge3d^2+9h-6\sqrt2\,|d|h^{1/2}=d^2+\big(\sqrt2\,|d|-3h^{1/2}\big)^2\ge d^2 .
\]
\end{proof}

\begin{lemma}\label{lem:energy}
Let $w\in A_2^d$, $c_I=b_I|I|^{-1/2}$ and $M=\|S_db\|_\infty$. Whenever $R\in\D$,
\begin{align*}
&\textup{(i)}\ \ \sum_{I\subseteq R}c_I^2\,\nu_I\le M^2\nu_R, &&\textup{(ii)}\ \ \sum_{I\subseteq R}d_I^2\,q_I\,\nu_I\le21\Aw\,\nu_R,\\
&\textup{(iii)}\ \ \sum_{I\subseteq R}d_I^2\,\sigma(I)\le\Aw\,\sigma(R), &&\textup{(iv)}\ \ \sum_{I\subseteq R}d_I^2\,w(I)\le4[w]_{A_\infty}w(R).
\end{align*}
\end{lemma}

\begin{proof}
Write $Q=\Aw$. Because every term is nonnegative, we need only bound, uniformly in $N$, the partial sums over $I\subseteq R$ satisfying $|I|>2^{-N}|R|$. The intervals $I\subseteq R$ satisfying $|I|=2^{-N}|R|$ will be called leaves.

(i) The Cauchy--Schwarz inequality yields $\nu_{I_l}+\nu_{I_r}\le\nu_I$. For leaves put $k_I=\nu_I$, and if $|I|>2^{-N}|R|$ put $k_I=\nu_I-\nu_{I_l}-\nu_{I_r}\ge0$. Thus, whenever $|J|\ge2^{-N}|R|$, $\nu_J=\sum_{K\subseteq J,\,|K|\ge2^{-N}|R|}k_K$ for $J\subseteq R$. Moreover, $\sum_{I\supseteq K}c_I^2\le S_db(x)^2$ for every $x\in K$, so
\[
\sum_{\substack{I\subseteq R\\|I|\ge2^{-N}|R|}}c_I^2\nu_I=\sum_{\substack{K\subseteq R\\|K|\ge2^{-N}|R|}}k_K\sum_{K\subseteq I\subseteq R}c_I^2\le M^2\sum_Kk_K=M^2\nu_R .
\]

(ii) Fix $I$. For $j\in\{l,r\}$, set $x_j=(w(I_j)/w(I))^{1/2}$ and $u_j=(\sigma(I_j)/\sigma(I))^{1/2}$. The vectors $x$ and $u$ are then unit vectors with nonnegative entries; moreover, $d_I=x_l^2-x_r^2$, $\nu_{I_j}=x_ju_j\nu_I$ and $q_{I_j}=4x_j^2u_j^2q_I$. Take $h=1-x\cdot u\ge0$ and $\Psi(J)=(21Q-q_J)\nu_J\ge0$. By Lemma~\ref{lem:twopoint} together with $q_I\le Q$ and $h\ge0$,
\[
\Psi(I)-\Psi(I_l)-\Psi(I_r)=\nu_I\Big[21Qh+q_I\Big(4\sum_j(x_ju_j)^3-1\Big)\Big]\ge q_I\nu_I\Big[21h+4\sum_j(x_ju_j)^3-1\Big]\ge d_I^2q_I\nu_I .
\]
Upon summing over $I\subseteq R$ satisfying $|I|>2^{-N}|R|$, telescoping bounds the left side by $\Psi(R)\le21Q\nu_R$.

(iii) As $\avg w_{I_l}=(1+d_I)\avg w_I$, $\avg w_{I_r}=(1-d_I)\avg w_I$ and $|I_l|=|I_r|=\frac12|I|$,
\[
\frac{|I_l|}{\avg w_{I_l}}+\frac{|I_r|}{\avg w_{I_r}}=\frac{|I|}{\avg w_I}\cdot\frac1{1-d_I^2}\ge(1+d_I^2)\frac{|I|}{\avg w_I}.
\]
Summing $d_I^2|I|/\avg w_I\le|I_l|/\avg w_{I_l}+|I_r|/\avg w_{I_r}-|I|/\avg w_I$ over those $I\subseteq R$ for which $|I|>2^{-N}|R|$ yields a quantity bounded above by the sum of $|L|/\avg w_L$ over the leaves $L$; moreover, $\avg w_L\avg\sigma_L\ge1$ implies $|L|/\avg w_L\le\sigma(L)$. In view of $\sigma(I)=q_I|I|/\avg w_I\le Q|I|/\avg w_I$, this establishes (iii).

(iv) Let $H(I)=\int_Iw\log(w/\avg w_I)$. Jensen's inequality ensures that this is nonnegative, while the bound below shows that it is finite. Taking $p=\frac12(1+d_I)$,
\[
H(I)-H(I_l)-H(I_r)=w(I)\,\phi(p),\qquad\phi(p)=p\log(2p)+(1-p)\log(2(1-p)).
\]
From $\phi(\frac12)=\phi'(\frac12)=0$ together with $\phi''(p)=1/(p(1-p))\ge4$ we obtain $\phi(p)\ge2(p-\frac12)^2=\frac12d_I^2$; telescoping then yields $\sum_{I\subseteq R,\,|I|>2^{-N}|R|}d_I^2w(I)\le2H(R)$.
 To estimate $H(R)$, take $m=\avg w_R$ and $\lambda\ge m$. Those maximal $K\in\D$ satisfying $K\subseteq R$ and $\avg w_K>\lambda$ are proper subintervals of $R$, whence $\avg w_K\le2\lambda$. They lie in $\{M^d(w\1_R)>\lambda\}$ and cover $R\cap\{w>\lambda\}$ except for a null set. Consequently, $|R\cap\{M^d(w\1_R)>\lambda\}|\ge(2\lambda)^{-1}\int_{R\cap\{w>\lambda\}}w$.
Applying $M^d(w\1_R)\ge m$ on $R$ and integrating over $\lambda\in(m,\infty)$,
\[
H(R)\le\int_Rw\log_+\frac wm\le2\int_R\big(M^d(w\1_R)-m\big)\le2[w]_{A_\infty}w(R)
\]
as a consequence of \eqref{eq:Ainf}. Thus (iv) is proved.
\end{proof}

\begin{proposition}\label{prop:34}
For each $b$ and each $w\in A_2^d$,
\[
\|\pi_b\|_{L^2(w)}\le C\,\Aw^{3/4}\,\|S_db\|_\infty ,
\]
where $C$ is an absolute constant.
\end{proposition}

\begin{proof}
By solidity and the finite-subfamily principle in Section~\ref{sec:prelim}, it suffices to consider $b$ having only finitely many nonzero coefficients.
Set $Q=\Aw$ and $M=\|S_db\|_\infty$. Proposition~\ref{prop:drift} yields $\pi_b=U_b+D_b$ with $\|U_b\|_{L^2(w)}\le2Q^{1/2}M$; in addition, whenever $f=\varphi\sigma$, we obtain $|D_bf|\le T(|\varphi|\sigma)$, where now
\[
T(\varphi\sigma)=\sum_{I\in\D}\lambda_I\avg{\varphi\sigma}_I\1_I,\qquad\lambda_I=|c_Id_I| .
\]
We check, for $T$ with $v=w$, the testing conditions in Theorem~\ref{thm:LSU}. Choose $R\in\D$ and set
\[
\beta_I=\lambda_Iq_I|I|=\lambda_I\avg\sigma_Iw(I)=\lambda_I\avg w_I\sigma(I),\qquad\Lambda_I=\sum_{J\subseteq I}\beta_J .
\]
Upon expanding the square, we group the pairs $J\subseteq I$ to obtain
\begin{equation}\label{eq:testid}
\int\Big(\sum_{I\subseteq R}\lambda_I\avg\sigma_I\1_I\Big)^2w=\sum_{I\subseteq R}\frac{2\beta_I\Lambda_I-\beta_I^2}{w(I)},\qquad
\int\Big(\sum_{I\subseteq R}\lambda_I\avg w_I\1_I\Big)^2\sigma=\sum_{I\subseteq R}\frac{2\beta_I\Lambda_I-\beta_I^2}{\sigma(I)},
\end{equation}
and, because $\beta_I\le\Lambda_I$, each right side is bounded below by $\sum_{I\subseteq R}\beta_I\Lambda_I/w(I)$ (respectively $/\sigma(I)$) and above by twice this sum.

Since $\beta_J=|c_J|\,|d_J|\,q_J^{1/2}\nu_J$, applying Lemma~\ref{lem:energy}(i),(ii) together with the Cauchy--Schwarz inequality yields $\Lambda_I\le(21Q)^{1/2}M\nu_I$. We also have $\beta_I\nu_I=|c_Id_I|\,q_I^{1/2}w(I)\sigma(I)$, hence
\[
\sum_{I\subseteq R}\frac{\beta_I\Lambda_I}{w(I)}\le\sqrt{21}\,QM\sum_{I\subseteq R}|c_Id_I|\,\sigma(I),\qquad
\sum_{I\subseteq R}\frac{\beta_I\Lambda_I}{\sigma(I)}\le\sqrt{21}\,QM\sum_{I\subseteq R}|c_Id_I|\,w(I).
\]
Applying the Cauchy--Schwarz inequality on each path, together with $\sum_{I\ni x}c_I^2\le M^2$ and Lemma~\ref{lem:energy}(iii), yields
\[
\sum_{I\subseteq R}|c_Id_I|\,\sigma(I)=\int_R\sum_{I\subseteq R,\,I\ni x}|c_Id_I|\,d\sigma(x)\le M\Big(\sigma(R)\sum_{I\subseteq R}d_I^2\sigma(I)\Big)^{1/2}\le Q^{1/2}M\sigma(R),
\]
and similarly Lemma~\ref{lem:energy}(iv) yields $\sum_{I\subseteq R}|c_Id_I|\,w(I)\le2[w]_{A_\infty}^{1/2}Mw(R)$. Combining \eqref{eq:testid} with \eqref{eq:Ainf} gives
\[
\begin{aligned}
\int\Big(\sum_{I\subseteq R}\lambda_I\avg\sigma_I\1_I\Big)^2w&\le2\sqrt{21}\,Q^{3/2}M^2\sigma(R),\\
\int\Big(\sum_{I\subseteq R}\lambda_I\avg w_I\1_I\Big)^2\sigma&\le4\sqrt{21}\,Q[w]_{A_\infty}^{1/2}M^2w(R)\le4\sqrt{21 c}\,Q^{3/2}M^2w(R).
\end{aligned}
\]
Thus $\mathfrak T+\mathfrak T^*\le\big((2\sqrt{21})^{1/2}+(4\sqrt{21})^{1/2}c^{1/4}\big)Q^{3/4}M\le C_1Q^{3/4}M$, while $\|D_b\|_{L^2(w)}\le C_1 C_TQ^{3/4}M$ follows from Theorem~\ref{thm:LSU}. Incorporating the bound for $U_b$ establishes the proposition with $C=2+C_1 C_T$.
\end{proof}

The exponent $\frac34$ arises from the first testing constant above, whereas the second is of order $Q^{1/2}[w]_{A_\infty}^{1/4}M$. For every symbol with bounded square function, Proposition~\ref{prop:34} improves on \eqref{eq:bez}; whenever $[w]_{A_\infty}>\Aw^{1/2}$, it also improves on Remark~\ref{rem:partial}.

\section{No universal rate}\label{sec:norate}

\begin{theorem}\label{thm:norate}
Given any $\rho:[2,\infty)\to(0,\infty)$ with $\rho(Q)/Q\to0$ as $Q\to\infty$, one can find $b\in\mathrm{CMO}(\R)\cap\cmo\cap L^\infty$ with $\limsup_{Q\to\infty}\Phi_b(Q)/\rho(Q)=\infty$. In particular, there exists a bounded $b\in\cmo$ belonging to no $\G_\alpha$ with $\alpha<1$.
\end{theorem}

\begin{proof}
Because $Q/\rho(Q)\to\infty$, choose $2\le Q_1<Q_2<\cdots$ satisfying $Q_n\to\infty$ and $Q_n/\rho(Q_n)\ge n^3$, and set $J_n=[10n+\frac14,10n+\frac12)$, a dyadic interval. On $J_n$, take $c_n$ to be the alternating chain from Proposition~\ref{prop:chain} with $\eps=1/n$, $\delta=1/Q_n$ and $L=\lceil Q_n\rceil$; denote the corresponding weight and test function by $w_n$ and $f_n=\sigma_n\1_{J_n}$.
We then have $[w_n]_{A_2^d}\le Q_n$, together with $\|\pi_{c_n}f_n\|_{L^2(w_n)}\ge Q_n\|f_n\|_{L^2(w_n)}/(\sqrt8\,n)$ by \eqref{eq:chainlower}.

\emph{Mollification.} Set $g_{n,j}=c_n*\phi_j$, with $\phi_j\ge0$ smooth, supported in $(-2^{-j},2^{-j})$, and having integral $1$, $j\ge3$.
We obtain $g_{n,j}\in C_c^\infty((10n+\frac18,10n+\frac58))$, $\int g_{n,j}=0$, $\|g_{n,j}\|_\infty\le\|c_n\|_\infty\le2/n$, and $g_{n,j}\to c_n$ in $L^1$; consequently, for every $I$, $(g_{n,j})_I\to(c_n)_I$.
For each fixed $n$, \eqref{eq:bez} yields uniform boundedness of the operators $\pi_{g_{n,j}}$ on $L^2(w_n)$, because $\|g_{n,j}\|_{\bmo}\lesssim1/n$.
When paired with the functions $w_n^{-1}h_I$ in $L^2(w_n)$, they give $(g_{n,j})_I\avg{f_n}_I\to(c_n)_I\avg{f_n}_I$.
The span of these functions is dense in $L^2(w_n)$. Indeed, orthogonality of $F\in L^2(w_n)$ to all of them implies $F\in L^1_{\mathrm{loc}}$, since for compact $K$ we have $\int_K|F|\le\|F\|_{L^2(w_n)}\sigma_n(K)^{1/2}$, and the Haar coefficients of $F$ vanish. Consequently, $F$ is constant on each half-line; moreover, $F=0$, since far out $w_n=1$.
It follows that $\pi_{g_{n,j}}f_n\rightharpoonup\pi_{c_n}f_n$ weakly in $L^2(w_n)$. Weak lower semicontinuity of the norm therefore gives some $g_n=g_{n,j(n)}$ satisfying $\|\pi_{g_n}f_n\|_{L^2(w_n)}\ge Q_n\|f_n\|_{L^2(w_n)}/(2\sqrt8\,n)$.

\emph{The symbol.} Set $b=\sum_ng_n$. Disjointness of the supports of the $g_n$, together with $\|g_n\|_\infty\le2/n$, implies that $b$ is bounded; its partial sums converge to $b$ in $L^\infty$ and therefore also in $\mathrm{BMO}$ and $\bmo$.
All partial sums belong to $C_c^\infty\subseteq\mathrm{CMO}(\R)$. Moreover, every $g_n$ belongs to $\cmo$. Indeed, its Haar coefficients obey $|(g_n)_I|\le\frac14|I|^{3/2}\|g_n'\|_\infty$; consequently, scales below $2^{-K}$ contribute no more than $O(2^{-K})$ to $\|g_n\|_{\bmo}$. Furthermore, $(g_n)_I=0$ holds for each dyadic $I$ satisfying $|I|\ge2$: the midpoint of any such interval is an integer, whereas $\operatorname{supp}g_n$ lies inside an open interval of length one bounded by consecutive integers, making $h_I$ constant on $\operatorname{supp}g_n$.
Thus, retaining the coefficients of $g_n$ at scales between $2^{-K}$ and $1$ gives a finite Haar sum that approximates $g_n$ in $\bmo$. Closedness of $\mathrm{CMO}(\R)$ and $\cmo$ now yields $b\in\mathrm{CMO}(\R)\cap\cmo$.

\emph{Localization.} Suppose $\avg{f_n}_I\ne0$; either $I\subseteq J_n$ or $I\supsetneq J_n$ follows. Under the first alternative, none of the other blocks meets $I$. Under the second, any interval intersecting both $J_n$ and another block's support has length no less than $9$, and thus no less than $16$. Its midpoint is therefore an integer, so, by the preceding reasoning, the coefficient of that other block vanishes. Consequently $\pi_bf_n=\pi_{g_n}f_n$ and $\Phi_b(Q_n)\ge Q_n/(2\sqrt8\,n)$, which gives $\Phi_b(Q_n)/\rho(Q_n)\ge n^2/(2\sqrt8)\to\infty$. To obtain the last claim, choose $\rho(Q)=Q/\log(e+Q)$.
\end{proof}

The jumps in the finite blocks $c_n$ exclude them from $\mathrm{VMO}$ and prevent approximation by smooth functions in $\mathrm{BMO}$. Moreover, measuring an approximation error in $\mathrm{BMO}$ would amplify it by a factor $Q$. To avoid this, the argument employs weak convergence with the weight fixed.

\section{Compactness and essential norms}\label{sec:compact}

\begin{theorem}\label{thm:ess}
\textup{(i)} For every $b$ and every $w\in A_2^d$,
\begin{equation}\label{eq:ess}
\|\pi_b\|_{\mathrm{ess},L^2(w)}\le C_B\Aw\dist_{\bmo}(b,\cmo).
\end{equation}

\textup{(ii)} One can find $b$ satisfying $\dist_{\bmo}(b,\cmo)>0$ for which $\|\pi_b\|_{\mathrm{ess},L^2(w)}\ge\frac14\Aw\dist_{\bmo}(b,\cmo)$ holds for weights $w$ having arbitrarily large $\Aw$.

\textup{(iii)} One can find $b$ satisfying $\dist_{\bmo}(b,\cmo)=1$ for which $\pi_b$ is compact on no $L^2(w)$, $w\in A_2^d$, while $\Phi_b(Q)=Q^{1/2}$ holds for every $Q\ge1$. In particular, $\G_{1/2}\not\subseteq\cmo$.
\end{theorem}

\begin{proof}
(i) Let $\mathcal F_K=\{I:2^{-K}\le|I|\le2^K,\ I\cap[-2^K,2^K]\ne\emptyset\}$ and $\mathcal E_K(b)=\|b-b_{\mathcal F_K}\|_{\bmo}$. The increase of $\mathcal F_K$ with $K$ implies that $\mathcal E_K(b)$ decreases. Because $b_{\mathcal F_K}$ is a finite Haar sum, we obtain $\dist_{\bmo}(b,\cmo)\le\mathcal E_K(b)$. Conversely, the coefficients of a finite Haar sum $g$ belong to $\mathcal F_K$ when $K$ is large; in that case, $b-b_{\mathcal F_K}$ restricts $b-g$, which yields $\mathcal E_K(b)\le\|b-g\|_{\bmo}$. Consequently, $\mathcal E_K(b)\downarrow\dist_{\bmo}(b,\cmo)$. The finite rank of $\pi_{b_{\mathcal F_K}}$ gives $\|\pi_b\|_{\mathrm{ess}}\le\|\pi_{b-b_{\mathcal F_K}}\|\le C_B\Aw\mathcal E_K(b)$.

(ii) Let $\beta_n=\sum_{k<n}(-1)^k|I_k|^{1/2}h_{I_k}$ with $I_k=[0,2^{-k})$, and $b=\eps\sum_{n\ge1}\beta_n(\cdot-10n)$. Integer translations preserve $\D$; moreover, the argument of Proposition~\ref{prop:rate-sharp} yields $\|b\|_{\bmo}\le\sqrt2\eps$. 
Every finite Haar sum omits the coefficient  $\eps|I|^{1/2}$ of $b$ at $I=[10n,10n+1)$ for some $n$ depending on that sum. Thus $\dist_{\bmo}(b,\cmo)\ge\eps$, and consequently $\eps\le\dist_{\bmo}(b,\cmo)\le\sqrt2\eps$. Given $Q\ge2$, choose $\delta=1/Q$ and equip $w$ with the profile from Proposition~\ref{prop:chain} on each block $P_n=[10n,10n+1)$, $n\in\Z$: take $w(x)=(x-10n)^{1-\delta}$ on $P_n$, with $w=1$ elsewhere. Whenever a dyadic interval $J$ satisfies $|J|<1$, it is contained in a unit interval $[j,j+1)$. Therefore, either $w\equiv1$ on $J$ or $J\subseteq P_n$ for some $n$; the proof of Proposition~\ref{prop:chain} then gives $\avg w_J\avg\sigma_J\le1/\delta$. Whenever a dyadic interval $J$ satisfies $|J|\ge1$, it is a union of unit intervals and hence either contains or misses each block. As $w(P_n)=1/(2-\delta)\le1$, this yields $\avg w_J\le1$. In the case $|J|=1$, we obtain $\avg w_J\avg\sigma_J\le1/(\delta(2-\delta))\le1/\delta$. In the case $|J|\ge2$, at most $\lceil|J|/10\rceil\le|J|/2$ blocks are contained in $J$, each having $\sigma$-mass $1/\delta$; therefore $\avg\sigma_J\le1+\frac1{2\delta}\le\frac1\delta$. We conclude that $\Aw\le Q$. For $n\ge Q$, the functions $f_n=\sigma\1_{P_n}$ have disjoint supports. All coefficients of $b$ belong to the blocks $P_m$, $m\ge1$. Consequently, no coefficients occur on the ancestors of $P_n$, and $\pi_bf_n$ agrees with the paraproduct of the $n$-th block, which is a chain of length $n\ge Q$. On $P_n$, $w$ agrees with the weight $w_\delta$ from Proposition~\ref{prop:chain} for $P=P_n$; thus \eqref{eq:chainlower} yields $\|\pi_bf_n\|\ge\eps Q\|f_n\|/\sqrt8$. Finally, $f_n/\|f_n\|\rightharpoonup0$ implies $\|\pi_b\|_{\mathrm{ess}}\ge\eps Q/\sqrt8\ge Q\dist_{\bmo}(b,\cmo)/4$.

(iii) Let $b=\sum_{j\in\Z}h_{[j,j+1)}$. Counting the unit intervals in $J$ gives the Carleson sum over $J$; hence $\|b\|_{\bmo}=1$, while the argument of (ii) yields $\dist_{\bmo}(b,\cmo)=1$. Lemma~\ref{lem:single} gives $\|\pi_b\|_{L^2(w)}=\sup_j(\avg w_{[j,j+1)}\avg\sigma_{[j,j+1)})^{1/2}\le\Aw^{1/2}$. Moreover, each of the infinitely many singular values is at least $1$, which shows that $\pi_b$ is not compact. Given $w_t=t$ on $[j,j+\frac12)$ and $w_t=1$ on $[j+\frac12,j+1)$, $j\in\Z$, we obtain $[w_t]_{A_2^d}=\frac{(t+1)^2}{4t}$ together with $\|\pi_b\|_{L^2(w_t)}=[w_t]_{A_2^d}^{1/2}$. As $\frac{(t+1)^2}{4t}$ attains every value $Q\ge1$, it follows that $\Phi_b(Q)=Q^{1/2}$.
\end{proof}

The approximation numbers are bounded by the truncation from (i) as well, $a_n(T)=\inf\{\|T-R\|:\operatorname{rank}R<n\}$.

\begin{proposition}\label{prop:approx}
Given $K\ge0$ and $w\in A_2^d$,
\begin{equation}\label{eq:approx}
a_{4^{K+2}+1}(\pi_b;L^2(w))\le C_B\Aw\|b^{>K}\|_{\bmo}+C_W'\sqrt{\min(2K+1,\Aw)\,\Aw}\;\mathcal E_K(b).
\end{equation}
\end{proposition}

\begin{proof}
For $|j|\le K$, the number of dyadic intervals of length $2^j$ meeting $[-2^K,2^K]$ is at most $2^{K+1-j}+2$. Consequently, $\#\mathcal F_K\le\sum_{j=-K}^K2^{K+1-j}+2(2K+1)<4^{K+1}+4K+2\le4^{K+2}$, with the rank of $\pi_{b_{\mathcal F_K}}$ at most $4^{K+2}$. Decompose $b-b_{\mathcal F_K}$ into the portions at scales $|k|>K$ and $|k|\le K$. Because $\mathcal F_K$ contains only scales $|k|\le K$, the first portion constitutes the scale tail of $b$; therefore \eqref{eq:bez} bounds its paraproduct by $C_B\Aw\|b^{>K}\|_{\bmo}$. For the second portion, the support lies on $2K+1$ consecutive scales, and its $\bmo$ norm is at most $\mathcal E_K(b)$. Applying Theorem~\ref{thm:window} with $[w]_{A_\infty},[\sigma]_{A_\infty}\le c\Aw$ therefore bounds its paraproduct by $C_W'\sqrt{\min(2K+1,\Aw)\Aw}\,\mathcal E_K(b)$.
\end{proof}

Thus \eqref{eq:approx} shows that a rate $\mathcal E_K(b)\lesssim\omega(K)$ gives a decay rate for the approximation numbers. For branching trees, however, the reverse comparison with best Haar truncations fails.

\begin{example}\label{ex:branch}
Let $b=h_{[0,1)}+2^{-1/2}(h_{[0,\frac12)}+h_{[\frac12,1)})$. Because $\avg f_{[0,1)}=\frac12(\avg f_{[0,\frac12)}+\avg f_{[\frac12,1)})$,
\[
\pi_bf=\avg f_{[0,\frac12)}\Big(\tfrac12h_{[0,1)}+2^{-\frac12}h_{[0,\frac12)}\Big)+\avg f_{[\frac12,1)}\Big(\tfrac12h_{[0,1)}+2^{-\frac12}h_{[\frac12,1)}\Big)
\]
has rank two, giving $a_3(\pi_b;L^2(w))=0$ for every weight. In contrast, if $F$ consists of at most two of the three intervals, then $b-b_F$ retains at least one of the coefficients of $b$. Each single term, $h_{[0,1)}$ or $2^{-1/2}h_I$ with $|I|=1/2$, has $\bmo$ norm; hence $\|b-b_F\|_{BMO_d}\geq 1$ for every such $F$.

\end{example}

\section{Concluding remarks}\label{sec:concl}

We investigated, for a fixed symbol, the growth of the weighted norm of the dyadic paraproduct in terms of the $A_2^d$ characteristic, as measured by $\Phi_b(Q)$ and the classes $\G_\alpha$. Every nonzero symbol supported on one scale attains the lower growth $Q^{1/2}$, whereas the upper growth $Q$ is attained over $\bmo$. Compactness does not determine where the growth falls between these bounds: although compactness of $\pi_b$ implies $\Phi_b(Q)=o(Q)$, it imposes no rate (Theorem~\ref{thm:norate}); moreover, minimal growth is possible for a non-compact symbol (Theorem~\ref{thm:ess}). We instead found two properties of the Haar coefficients that govern growth and determine it for the model cases considered here. One concerns their distribution over scales. For symbols supported on $L$ consecutive scales, the growth is like $\sqrt{\min(L,Q)Q}$; polynomial decay of the scale tails yields the sharp exponent $\max(1-\kappa,\frac12)$ for $\kappa\ne\frac12$ (Theorem~\ref{thmA}). At $\kappa=\frac12$, however, the worst growth is between $(Q\log Q)^{1/2}$ and $Q^{1/2}\log Q$ (Theorem~\ref{thmA} and Corollary~\ref{cor:critical}). The other property concerns nesting of the Haar support. On a chain, the relevant quantity is the $\ell^2$ sum of the normalized coefficients. When the coefficients are nonincreasing, growth is comparable to $Q^{1/2}$ times the $\ell^2$ sum of their first $Q$ terms (Theorem~\ref{thm:monotone}). For every symbol, boundedness of $S_db$ is necessary for minimal growth (Theorem~\ref{thmB}), while it is sufficient for the sublinear growth $Q^{3/4}$ (Proposition~\ref{prop:34}).

These results leave the description of general branching Haar supports incomplete. Whether the converse of Theorem~\ref{thmB}(ii) is valid remains the principal open question: does \eqref{eq:Sq} hold for every symbol and every weight? The reasons presented in Section~\ref{sec:Sq} lead us to expect an affirmative answer.

\begin{conjecture}\label{conj:S}
$\|\pi_b\|_{L^2(w)}\le C\Aw^{1/2}\|S_db\|_\infty$ for every $b$ and every $w\in A_2^d$. By Theorem~\ref{thmB}, this is equivalent to $\G_{1/2}=\{b:S_db\in L^\infty\}$ with $[b]_{\G_{1/2}}\simeq\|S_db\|_\infty$.
\end{conjecture}

What we establish is as follows. For chain symbols, the conjectured bound is valid at every weight (Theorem~\ref{thmB}). It is also valid for every symbol over weights with $A_2$ imbalance of bounded energy off a single chain or off a family of chains with disjoint tops (Theorems~\ref{thm:spine} and~\ref{thm:forest}). Every weight employed for the lower bounds in this paper belongs to this class (Corollary~\ref{cor:power}, Remark~\ref{rem:blocks}). In the case of general weights, precisely one extra factor $[w]_{A_\infty}^{1/2}$ remains in the argument of Remark~\ref{rem:partial}; it arises from the lower square function estimate. The growth $Q^{3/4}$ is given by Proposition~\ref{prop:34}. Consequently
\[
\G_{1/2}\subseteq\{b:S_db\in L^\infty\}\subseteq\G_{3/4},
\]
where strictness of the second inclusion follows from $b_{1/2}\in\G_{3/4}$ and $S_db_{1/2}\notin L^\infty$, by Corollary~\ref{cor:critical}. The first inclusion, according to the conjecture, is an equality. Remark~\ref{rem:partial} implies that any counterexample must involve a sequence of weights for which $[w]_{A_\infty}$ is unbounded. Moreover, Theorem~\ref{thm:forest} requires this sequence to escape every class of weights having uniformly bounded imbalance energy off some family of chains with disjoint tops. One possible mechanism consists of singular chains nested across scales, although other mechanisms are not ruled out by the present results. An affirmative answer would give a description of $\G_{1/2}$ independent of weights. It would also permit replacing the exponent $\min(\frac32-\alpha,\frac34)$ from Remark~\ref{rem:rates} with $\frac12$, because $\|S_db\|_\infty^2\le\sum_k\|\Delta_kb\|_\infty^2$, and resolve the critical case of Theorem~\ref{thmA} in the following way.

\begin{remark}\label{rem:critical}
Suppose that Conjecture~\ref{conj:S} is valid, with $\|b^{>K}\|_{\bmo}\le C_0(1+K)^{-1/2}$ for every $K\ge0$. Take $N=\lceil Q\rceil$ and $b^{\le N}=b-b^{>N}$. A given point belongs to exactly one interval at each scale; Lemma~\ref{lem:basic}(i) therefore yields
\[
\|S_db^{\le N}\|_\infty^2\le\sum_{|k|\le N}\|\Delta_kb\|_\infty^2\le\|b\|_{\bmo}^2+2C_0^2\sum_{k=1}^Nk^{-1}\le\|b\|_{\bmo}^2+2C_0^2(1+\log N),
\]
while the paraproduct of $b^{>N}$ on $L^2(w)$, $\Aw\le Q$, is bounded by $C_BC_0Q(1+N)^{-1/2}\le C_BC_0Q^{1/2}$ using \eqref{eq:bez}. Thus $\Phi_b(Q)\le C(\|b\|_{\bmo}+C_0)(Q\log(e+Q))^{1/2}$. Corollary~\ref{cor:critical} shows that removing the factor $(\log Q)^{1/2}$ is impossible. Accordingly, the conjecture would specify, up to constants, the growth at the critical rate.
\end{remark}

The point at which the factor $Q^{1/4}$ between the result and the conjecture is lost can be seen in the proof of Proposition~\ref{prop:34}. Already, the orthogonal part $U_b$ satisfies the conjectured bound. Meanwhile, the positive operator $T$ dominates the drift, and Theorem~\ref{thm:LSU} makes its norm comparable to its testing constants. Taking $\beta_I$ and $\Lambda_I$ as in the proof of Proposition~\ref{prop:34} and $M=\|S_db\|_\infty$, one sees from \eqref{eq:testid} that the following inequalities would imply Conjecture~\ref{conj:S}
\begin{equation}\label{eq:sharptest}
\sum_{I\subseteq R}\frac{\beta_I\Lambda_I}{w(I)}\le C\Aw M^2\sigma(R),\qquad\sum_{I\subseteq R}\frac{\beta_I\Lambda_I}{\sigma(I)}\le C\Aw M^2w(R),\qquad R\in\D,
\end{equation}
which depend solely on the normalized coefficients and on the weight’s dyadic statistics. The signs of both $c_I$ and $d_I(w)$ are disregarded by $T$; consequently, \eqref{eq:sharptest} gives a sufficient condition stronger than the conjecture. For general branching supports, its validity remains unestablished here, although it agrees with every case resolved in this paper.
For chain symbols, it holds at every weight with an absolute constant: the proof of Theorem~\ref{thm:chain} relies solely on the pointwise bound \eqref{eq:ptchain}, and $T$ satisfies that bound with $\lambda_I\le|c_I|$ replacing $|c_I|$. Moreover, combining this observation for the chain part with the testing estimates in the proofs of Theorems~\ref{thm:spine} and~\ref{thm:forest} establishes it for every symbol at the weights of those theorems, with a constant depending on $E$.

In proving Proposition~\ref{prop:34}, the replacement for \eqref{eq:sharptest} is a product arising from two applications of the Cauchy--Schwarz inequality. The application across generations below $I$ bounds $\Lambda_I$ by Lemma~\ref{lem:energy}(i),(ii); the application along paths uses Lemma~\ref{lem:energy}(iii),(iv). A factor $Q^{1/2}$ may be lost in either application, although the situations producing these losses are opposite.
Consider the weight $w_\delta$ from Proposition~\ref{prop:chain} on a chain with $\delta=1/Q$, so that along this chain $q_I\simeq Q$ and $|d_I|\simeq1$ hold.
The bound for $\Lambda_I$ is attained by a single nonzero coefficient, whereas the path inequality incurs a loss of order at least $Q^{1/2}$. Conversely, taking coefficients of size $Q^{-1/2}$ on the first $Q$ generations attains the path inequality but incurs a loss of order $Q^{1/2}$ in the bound for $\Lambda_I$.
For both choices, the two testing sums have orders $QM^2\sigma(R)$ and $QM^2w(R)$. This is required by the chain case of \eqref{eq:sharptest} because $[w_\delta]_{A_2^d}\le Q$. These observations suggest handling the two steps jointly, for example by a stopping time adapted both to the levels of $q_I$ and to how the mass $\sum_Ic_I^2$ is distributed along paths.
We intend to investigate this approach alongside refinements of the decomposition in Proposition~\ref{prop:drift} using stopping-time arguments adapted simultaneously to $w$ and $\sigma$.

Further avenues concern lower bounds for the approximation numbers of $\pi_b$ on $L^2(w)$ incorporating the redundancy of branching trees, along the lines of the Schatten class results in \cite{LLW24}. Other directions are two-weight and Bloom versions \cite{LL22}, matrix weights \cite{LLSW26}, and commutators $[b,T]$ with Calder\'on--Zygmund operators. For these commutators, the linear dependence is replaced by the quadratic dependence of \cite{Ch11,CPP12}; the corresponding questions concerning scale profile and nesting will be addressed separately.

\section*{Declarations}

\noindent\textbf{Funding.} This work received support through the National Research Foundation of Korea (NRF) grant funded by the Korea government (MSIT) (RS-2026-25602221).

\smallskip\noindent\textbf{Competing interests.} The author declares no competing interests relevant to this article's content.

\smallskip\noindent\textbf{Data availability.} Data sharing does not apply to this article because the current study neither generated nor analysed any datasets.

\smallskip\noindent\textbf{Use of AI tools.} During the preperation of this manuscripts, the author used generative AI(Claude Opus 5.5) solely for the purpose of language polishing and smoothing the English text. After usinf this tool, the author carefully reviewed and edited the content as needed ans assumes full reponsibility for the funal publication.

\end{document}